\documentclass[11pt,leqno]{article}
\usepackage{graphicx, amsfonts, amsthm, amsxtra, amssymb, verbatim,latexsym,bbm}
\usepackage{addfont}
\usepackage{hyperref}
\usepackage[mathscr]{euscript}
\usepackage[all,cmtip]{xy}

\begin{document}

\def\Z{\mathbb{Z}}                   
\def\Q{\mathbb{Q}}                   
\def\C{\mathbb{C}}                   
\def\N{\mathbb{N}}                   
\def\dR{{\rm dR}}                    
\def\Gm{\mathbb{G}_m}                 
\def\Ga{\mathbb{G}_a}                 
\def\Tr{{\rm Tr}}                      
\def\tr{{{\mathsf t}{\mathsf r}}}                 
\def\spec{{\rm Spec}}            
\def\ker{{\rm ker}}              
\def\GL{{\rm GL}}                
\def\k{{\sf k}}                     
\def\ring{{ R}}                   
\def\X{{\sf X}}                      
\def\T{{\sf T}}                      
\def\Ts{{\sf S}}
\def\cmv{{\sf M}}                    
\def\BG{{\sf G}}                       
\def\podu{{\sf pd}}                   
\def\ped{{\sf U}}                    
\def\per{{\sf  P}}                   
\def\gm{{\sf  A}}                    
\def\gma{{\sf  B}}                   
\def\ben{{\sf b}}                    

\def\Rav{{\mathfrak M }}                     
\def\Ram{{\mathfrak C}}                     
\def\Rap{{\mathfrak G}}                     

\def\nov{{\sf  n}}                    
\def\mov{{\sf  m}}                    
\def\Yuk{{\sf Y}}                     
\def\Ra{{\sf R}}                      
\def\hn{{ h}}                      
\def\cpe{{\sf C}}                     
\def\g{{\sf g}}                       
\def\t{{\sf t}}                       
\def\pedo{{\sf  \Pi}}                  

\def\Der{{\rm Der}}                   
\def\MMF{{\sf MF}}                    
\def\codim{{\rm codim}}                
\def\dim{{\rm    dim}}                
\def\Lie{{\rm Lie}}                   
\def\gg{{\mathfrak g}}                

\def\u{{\sf u}}                       

\def\imh{{  \Psi}}                 
\def\imc{{  \Phi }}                  
\def\stab{{\rm Stab }}               
\def\Vec{{\rm Vec}}                 
\def\prim{{\rm prim}}                  

\def\Fg{{\sf F}}     
\def\hol{{\rm hol}}  
\def\non{{\rm non}}  
\def\alg{{\rm alg}}  

\def\bcov{{\rm \O_\T}}       

\def\leaves{{\cal L}}        

\def\GM{{\rm GM}}

\def\perr{{\sf q}}        
\def\perdo{{\cal K}}   
\def\sfl{{\mathrm F}} 
\def\sp{{\mathbb S}}  

\newcommand\diff[1]{\frac{d #1}{dz}} 
\def\End{{\rm End}}              

\def\sing{{\rm Sing}}            
\def\cha{{\rm char}}             
\def\Gal{{\rm Gal}}              
\def\jacob{{\rm jacob}}          
\def\tjurina{{\rm tjurina}}      
\newcommand\Pn[1]{\mathbb{P}^{#1}}   
\def\Ff{\mathbb{F}}                  

\def\O{{\cal O}}                     
\def\as{\mathbb{U}}                  
\def\ring{{ R}}                         
\def\R{\mathbb{R}}                   

\def\Mat{{\rm Mat}}              
\def\cl{{\rm cl}}                

\def\hc{{\mathsf H}}                 
\def\Hb{{\cal H}}                    
\def\pese{{\sf P}}                  

\def\PP{\tilde{\cal P}}              
\def\K{{\mathbb K}}                  

\def\M{{\cal M}}
\def\RR{{\cal R}}
\newcommand\Hi[1]{\mathbb{P}^{#1}_\infty}
\def\pt{\mathbb{C}[t]}               
\def\W{{\cal W}}                     
\def\gr{{\rm Gr}}                
\def\Im{{\rm Im}}                
\def\Re{{\rm Re}}                
\def\depth{{\rm depth}}
\newcommand\SL[2]{{\rm SL}(#1, #2)}    
\newcommand\PSL[2]{{\rm PSL}(#1, #2)}  
\def\Resi{{\rm Resi}}              

\def\L{{\cal L}}                     
\def\Aut{{\rm Aut}}              
\def\any{R}                          
\newcommand\ovl[1]{\overline{#1}}    

\newcommand\mf[2]{{M}^{#1}_{#2}}     
\newcommand\mfn[2]{{\tilde M}^{#1}_{#2}}     

\newcommand\bn[2]{\binom{#1}{#2}}    
\def\ja{{\rm j}}                 
\def\Sc{\mathsf{S}}                  
\newcommand\es[1]{g_{#1}}            
\newcommand\V{{\mathsf V}}           
\newcommand\WW{{\mathsf W}}          
\newcommand\Ss{{\cal O}}             
\def\rank{{\rm rank}}                
\def\Dif{{\cal D}}                   
\def\gcd{{\rm gcd}}                  
\def\zedi{{\rm ZD}}                  
\def\BM{{\mathsf H}}                 
\def\plf{{\sf pl}}                             
\def\sgn{{\rm sgn}}                      
\def\diag{{\rm diag}}                   
\def\hodge{{\rm Hodge}}
\def\HF{{\sf F}}                                
\def\WF{{\sf W}}                               
\def\HV{{\sf HV}}                                
\def\pol{{\rm pole}}                               
\def\bafi{{\sf r}}
\def\id{{\rm id}}                               
\def\gms{{\sf M}}                           
\def\Iso{{\rm Iso}}                           

\def\hl{{\rm L}}    
\def\imF{{\rm F}}
\def\imG{{\rm G}}
\def\cy{{Calabi-Yau }}
\def\DHR{{\rm DHR }}            
\def\H{{\sf{ E}}}
\def\P {\mathbb{P}}                  
\def\E{{\cal E}}
\def\L{{\sf L}}             
\def\CX{{\cal X}}
\def\dt{{\sf d}}             
\def\LG{{\sf G}}   
\def\LA{{\rm Lie}(\LG)}   
\def\amsy{\mathfrak{G}}  
\def\gG{{\sf g}}   
\def\gL{\mathfrak{g}}   
\def\rvf{{\sf H}}   
\def\cvf{{\sf F}}   
\def\di2{{ m}}   
\def\a{\mathfrak{a}}
\def\b{\mathfrak{b}}

\def\mmat{e}
\def\cmat{f}
\def\rmat{h}

\def\hn{{\sf h}}                      

\def\tgtwo{\mathfrak{t}}
\def\tgone{\mathfrak{t}}
\def\tgzero{\mathfrak{t}}
\def\kgone{\mathfrak{k}}
\def\ggtwo{\mathfrak{g}}
\def\ggzero{\mathfrak{g}}

\def\Yukk{{\sf C}^\alg}

\def\sl2{\mathfrak{sl}_2(\C)}

\def\cyqmfs{{ \widetilde{\mathcal{M}}}}             
\def\cymfs{{\mathcal{M}}}                            
\def\Ramvf{{\sf Ra}}             

\newtheorem{theo}{Theorem}[section]
\newtheorem{exam}{Example}[section]
\newtheorem{coro}{Corollary}[section]
\newtheorem{defi}{Definition}[section]
\newtheorem{prob}{Problem}[section]
\newtheorem{lemm}{Lemma}[section]
\newtheorem{prop}{Proposition}[section]
\newtheorem{rem}{Remark}[section]
\newtheorem{obs}{Observation}[section]
\newtheorem{conj}{Conjecture}
\newtheorem{nota}{Notation}[section]
\newtheorem{ass}{Assumption}[section]
\newtheorem{calc}{}
\numberwithin{equation}{section}

\begin{center}
{\LARGE\bf  Modular vector fields attached to Dwork family: $\mathfrak{sl}_2(\C)$ Lie algebra }
\footnote{ MSC2010:
32M25,  	
37F99,  	    
14J15,      
14J32.      
\\
Keywords: Complex vector fields, Gauss-Manin connection, Dwork family, Hodge filtration, modular form. }
\\

\vspace{.25in} {\large {\sc 
Younes Nikdelan \footnote{Universidade do Estado do
Rio de Janeiro (UERJ), Instituto de Matem\'atica
e Estat\'{i}stica (IME), Departamento de An\'alise Matem\'atica: Rua São Francisco Xavier, 524, Rio de Janeiro, Brazil / CEP: 20550-900. e-mail: younes.nikdelan@ime.uerj.br}}} \\
\end{center}
\vspace{.25in}
\begin{abstract}
This paper aims to show that a certain moduli space $\T$, which arises from the so-called
Dwork family of Calabi-Yau $n$-folds, carries
a special complex Lie algebra containing a copy of
$\sl2$. In order to achieve this goal,
we introduce an algebraic group $\LG$ acting from the right on $\T$ and describe its Lie
algebra $\LA$. We observe that $\LA$ is isomorphic to a Lie subalgebra of the space
of the vector fields on $\T$. In this way,
it turns out that $\LA$  and the modular
vector field $\Ra$ generate another Lie algebra
$\amsy$, called AMSY-Lie algebra, satisfying $\dim (\amsy)=\dim (\T)$. We find a
copy of $\mathfrak{sl}_2(\C)$ containing $\Ra$ as a Lie subalgebra of
$\amsy$. The proofs are based on an algebraic
method calling "Gauss-Manin connection in disguise". Some explicit examples for $n=1,2,3,4$ are stated as well.  \\
\end{abstract}

\section{Introduction}\label{section introduction}

A modular vector field $\Ra$ is a special (and unique) vector field existing on a certain
quasi-affine variety $\T$ satisfying a definite equation involving the
Gauss-Manin connection. In the present work, the quasi-affine variety $\textsf{T}$ is the moduli space
of the pairs formed by certain Calabi-Yau $n$-folds (arising from the Dwork
family) along with $n+1$ differential $n$-forms. The modular vector field $\Ra$, in some sense, can be considered as a
generalization of the systems of differential equations introduced by  G.
Darboux \cite{da78}, G. H. Halphen \cite{ha81} and S. Ramanujan
\cite{ra16}, which have been discussed in \cite{ho14}. A. Guillot, J. Rebelo and H. Reis in \cite{gui07},
\cite{guireb} and \cite{rr14} treat the semi-completeness and the
dynamics of the Halphen type vector field $H$, which is a generalization of the Darboux-Halphen vector field,
in a 3-dimensional complex variety $W$. They use the action of the Lie
group $\rm SL(2,\C)$ on $W$ and its Lie algebra
$\mathfrak{sl}_2(\C)$ to reduce the dynamics of $H$, and then employ
a phenomenon of the integration which allows one to understand
the dynamics in detail. More specifically, $H$ together with a
constant vector field and the radial vector field generates a Lie algebra isomorphic to
$\mathfrak{sl}_2(\C)$.

In this article we describe an algebraic group action on $\T$,  its Lie algebra and, in
particular,  a copy of $\mathfrak{sl}_2(\C)$ containing $\Ra$. These should provide
foundations for the dynamical study of $\Ra$ along with the nature of its solutions.

From the early days of \textit{Mirror Symmetry} ideas revolving around \textit{Calabi-Yau manifolds}, a
vast number of works in mathematics and physics have been devoted to these
subjects. A pioneer work on the development of these theories is the 1991's paper of Candelas et al. \cite{can91} where they
used the mirror symmetry to predict the number of rational curves on
quintic 3-folds. These numbers come from the
coefficients of the $q$-expansion of a function called \emph{Yukawa
coupling}. H. Movasati in \cite{ho22} reencountered the Yukawa coupling
by applying an algebraic method, called \emph{Gauss-Manin connection
in disguise}, GMCD for short. GMCD simplifies the Candelas et al.'s method since
it does not require to compute either the periods or the variation of
the Hodge structure: all these objects become encoded in
the Gauss-Manin connection. H. Movasati started the project GMCD by applying to elliptic curves \cite{ho14}.  This
project, so far, contains a sequence of interesting completed works, see for instance
\cite{GMCD-MQCY3}, \cite{GMCD-NLH}, \cite{alimov}, \cite{movnik}.
One of the main objects of GMCD is the existence of a certain (unique) vector
field $\Ra$ on the moduli space  $\T$ of a special family of \cy
varieties enhanced with the differential forms. In lower dimensions, where so far the
solutions of $\Ra$ have been determined, it turns out that the
$q$-expansion of the solutions, up to multiplying by a constant rational number, has
integer coefficients. In particular, in the
case of elliptic curves and $K3$ surfaces the solutions
can be written in terms of (quasi-)modular forms (see
\cite{ho14} and the cases $n=1,2$ of \cite{movnik}). The vector field $\Ra$ is called \emph{modular vector field} and the $q$-expansion
of its solutions are called \emph{Calabi-Yau modular forms} by H. Movasati. One of the objectives of the project
GMCD, considering the modularity of the \cy varieties, is  somehow to generalize the theory of modular forms. In this generalization the \cy modular forms
can be regarded as the suitable candidates to substitute the classical modular forms, see for instance \cite{GMCD-MQCY3}.

On one hand, the mirror quintic 3-fold is the main example of
\cite{alimov} in which the authors describe a special Lie Algebra on the
moduli space of the non-rigid compact Calabi-Yau threefolds enhanced with differential
forms and discuss its relation to the Bershadsky-Cecotti-Ooguri-Vafa
holomorphic anomaly equation (see also \cite{ali17}).  On the other hand, the mirror
quintic 3-fold is the particular case $n=3$ of the family of the $n$-dimensional
Calabi-Yau varieties $X=X_\psi,\ \psi\in \Pn 1-\{0,1,\infty\}$,
arising from the so-called Dwork family:
\begin{equation}
\label{12jan2016}
x_0^{n+2}+x_1^{n+2}+\ldots+x_{n+1}^{n+2}-(n+2)\psi
x_0x_1\ldots x_{n+1}=0.
\end{equation}
The enhanced moduli space $\T$ and the relevant modular vector field $\Ra$ of this family have been
discussed in \cite{movnik}. These facts lead us to extend the Lie algebras discussed in \cite{alimov} to the enhanced moduli space $\T$
arising from the Dwork family for any positive integer $n$.

From now on, unless otherwise is stated, by \emph{mirror variety} $X$ we mean the \cy
$n$-fold $X=X_\psi,\ \psi\in \Pn 1 \setminus\{0,1,\infty\}$, obtained by the desingularization
of the quotient space of the \cy varieties of the Dwork family \eqref{12jan2016} under a
group action (see \cite[\S2]{movnik}).  By \emph{enhanced moduli space} $\T=\T_n$ we mean the moduli space of
the pairs $(X,[\alpha_1,\alpha_2,\ldots,\alpha_{n+1}])$ where $X$ is the mirror variety and  $\{\alpha_1,\alpha_2,\ldots,\alpha_{n+1}\}$
refers to a basis of the $n$-th algebraic de Rham cohomology $H_\dR^n(X)$ satisfying the following two conditions:  (i) the basis must be compatible with the
Hodge filtration of $H_\dR^n(X)$, (ii)  the intersection form matrix of the basis has to be constant (see Section \ref{section MSGMCD}).
It is well known that $\dim( H^n_\dR(X))=n+1$ and the Hodge numbers $h^{ij},\ i+j=n,$ of $X$ are all one.  The universal family $\pi:\X\to\T$
together with global sections $\alpha_i,\ \ i=1,\cdots, n+1$, of the \emph{relative algebraic de
Rham cohomology} $H^n_\dR(\X/\T)$ were constructed in \cite{movnik}. In the same reference, $\T$ was endowed with a
complete chart $\t=(t_1,t_2,\ldots,t_\dt)$ where $\dt=\dt_n:=\dim ( \T)$ (see \eqref{29dec2015}).
Let
$$
\nabla:H_{\dR}^{n}(\X/\T)\to
\Omega_\T^1\otimes_{\O_\T}H_{\dR}^{n}(\X/\T),
$$
be the algebraic \emph{Gauss-Manin connection} on $H^n_\dR(\X/\T)$ where $\O_\T$ is the $\C$-algebra of the regular functions on
$\T$ and $\Omega_\T^1$ refers to the $\O_{\T}$-module of the differential 1-forms on $\T$. The reader who is not familiar with the Gauss-Manin connection can substitute $\nabla$ with $d\int_{\delta_{\t}}$ where  $\t=(t_1,t_2,\ldots,t_{\dt})$ and $\delta_\t$ is an $n$-dimensional homology class in $\X_\t:=\pi^{-1}(\t)$. In fact, the equality $\nabla \alpha_i=\sum_{j=1}^{n+1}\zeta_{ij}\otimes\alpha_j$ is equivalent to saying that $d\int_{\delta_{\t}}\alpha_i=\sum_{j=1}^{n+1} \zeta_{ij}\int_{\delta_{\t}}\alpha_j$ where $\zeta_{ij}$'s are one forms.
For any $\H\in \mathfrak{X}(\T)$, in which $\mathfrak{X}(\T)$ stands for the Lie algebra of vector fields on $\T$, the linear application
$\nabla_{\H}:H_{\dR}^{n}(\X/\T)\to H_{\dR}^{n}(\X/\T)$ refers to the Gauss-Manin connection $\nabla$ composed with the map $\H\otimes \textrm{id}:\Omega_\T^1\otimes_{\O_\T}H_{\dR}^{n}(\X/\T)\to H_{\dR}^{n}(\X/\T)$. We
attach an $(n+1)\times(n+1)$ matrix $\gm_\H$, called
\emph{Gauss-Manin connection matrix}, to the vector field $\H$
defined as follow:
\begin{equation}\label{eq A_H}
\nabla_\H\alpha=\gm_\H \alpha,
\end{equation}
in which
\begin{equation} \label{eq alpha}
\alpha:=\left(
{\begin{array}{*{20}{c}} {{\alpha _1}}&{{\alpha _2}}& \ldots
&{{\alpha _{n + 1}}}
\end{array}} \right)^\tr.
\end{equation}
In the rest of this article the notation $\alpha$, being clear in the context, either is used for the columnar matrix \eqref{eq
alpha}, or for the $(n+1)$-tuple
$(\alpha_1,\alpha_2,\ldots,\alpha_{n+1})$.

Here, before continuing with the technical parts and the methodology, we state some explicit examples of Lie algebras containing $\Ra$ which are isomorphic to
$\sl2$. First, notice that the special linear Lie algebra $\sl2$ is the  Lie algebra of $2\times 2$ matrices with trace zero.
Three matrices
\begin{equation}\label{eq sbsl2}
\mmat:=\left(
        \begin{array}{cc}
          0 & 1 \\
          0 & 0 \\
        \end{array}
      \right)\ \ , \ \ \ \
\cmat:=\left(
        \begin{array}{cc}
          0 & 0 \\
          1 & 0 \\
        \end{array}
      \right)\ \ , \ \ \ \
\rmat:=\left(
        \begin{array}{cc}
          1 & 0 \\
          0 & -1 \\
        \end{array}
       \right)\ ,
\end{equation}
form the standard basis of $\sl2$ with the commutators:
\begin{equation}\label{eq lbsl2}
[\mmat,\cmat]=\rmat \ \ , \ \ \ \ [\rmat,\mmat]=2\mmat \ \ , \ \ \ \ [\rmat,\cmat]=-2\cmat\, .
\end{equation}
From now on, the notations $\mmat,\cmat,\rmat$ are saved for the standard basis \eqref{eq sbsl2} of $\sl2$
with the Lie brackets \eqref{eq lbsl2}. In order to prove that a determined Lie algebra is isomorphic to
$\sl2$ it is enough to introduce the corresponding standard basis. It is well known, see for instance \cite[\S2.12]{GMCD-MQCY3}, that
the Ramanujan vector field
\begin{equation}\label{eq ram vf}
\Ra=-(t_1^2-\frac{1}{12}t_2)\frac{\partial}{\partial t_1}-(4t_1t_2-6t_3)\frac{\partial}{\partial t_2}-(6t_1t_3-\frac{1}{3}t_2^2)\frac{\partial}{\partial t_3}\ ,
\end{equation}
together with the vector fields
\[
\rvf=2t_1\frac{\partial}{\partial t_1}+4t_2\frac{\partial}{\partial t_2}+6t_3\frac{\partial}{\partial t_3}\quad , \qquad \cvf=\frac{\partial}{\partial t_1}\ ,
\]
is isomorphic to $\sl2$ under the correspondences $\Ra\mapsto e, \ \cvf\mapsto f, \ \rvf\mapsto h$. In the point of view of GMCD, the vector fields $\Ra,\rvf,\cvf$ satisfy
{ \[
\nabla_\Ra\alpha=\left(
                    \begin{array}{cc}
                      0 & 1 \\
                      0 & 0 \\
                    \end{array}
                  \right) \alpha \ , \ \ \nabla_\rvf \alpha=\left(
                    \begin{array}{cc}
                      1 & 0 \\
                      0 & -1 \\
                    \end{array}
                  \right) \alpha \ , \ \ \nabla_\cvf\alpha=\left(
                    \begin{array}{cc}
                      0 & 0 \\
                      1 & 0 \\
                    \end{array}
                  \right) \alpha \ ,
\] }
where $\alpha=(\ \alpha_1\ \ \alpha_2\ )^\tr$ and $\nabla$ is the Gauss-Manin connection of the universal family
of the elliptic curves
\[
y^2=4(x-t_1)^3-t_2(x-t_1)-t_3 \ , \ \ \alpha_1=[\frac{dx}{y}], \ \alpha_2=[\frac{xdx}{y}], \ \text{with} \ \ 27t_3^2-t_2^3\neq 0\ .
\]
A similar statement holds for the Darboux-Halphen vector field, see \cite{ho12,gui07}. Analogously, on the enhanced moduli space $\T=\T_n$ arising from the Dwork family we find a unique modular vector field $\Ra\in \mathfrak{X}(\T)$ which along with two other vector fields
forms a copy of $\sl2$.
\begin{exam} \label{exam n=1,4}
In this example we explicitly state the modular vector field $\Ra$ and vector fields $\rvf$ and $\cvf$ on the enhanced moduli space $\T$ for $n=1,4$. These three vector fields form a copy of $\sl2$  under the correspondences $\Ra\mapsto e, \ \cvf\mapsto f, \ \rvf\mapsto h$. For details of the computations and more examples see Section \ref{subsection sl2} and Section \ref{subsection cy3}.
\begin{itemize}
  \item {\bf $n=1$}:
{\small
  \begin{align}
  &\Ra=(-t_1t_2-9(t_1^3-t_3))\frac{\partial}{\partial
t_1}+(81t_1(t_1^3-t_3)-t_2^2)\frac{\partial}{\partial
t_2}+(-3t_2t_3)\frac{\partial}{\partial t_3}, \label{eq mvf1}\\
&\rvf=t_1\frac{\partial}{\partial
t_1}+2t_2\frac{\partial}{\partial
t_2}+3t_3\frac{\partial}{\partial t_3},\label{eq rvf1}\\
&\cvf=\frac{\partial}{\partial t_2}\, , \label{eq cvf1}
  \end{align}}
\item $n=4$:
  {\small
\begin{align}
\Ra&=(t_3-t_1t_2)\frac{\partial}{\partial
t_1}+\frac{6^{-2}t_3^2t_4t_8-t_1^6t_2^2+t_2^2t_6}{t_1^6-t_6}\,\frac{\partial}{\partial
t_2} \label{eq mvf4} \\ &+\frac{6^{-2}t_3^2t_5t_8-3t_1^6t_2t_3+3t_2t_3t_6}{t_1^6-t_6}\,\frac{\partial}{\partial
t_3}
+\frac{-6^{-2}t_3^2t_7t_8-t_1^6t_2t_4+t_2t_4t_6}{t_1^6-t_6}\,\frac{\partial}{\partial
t_4} \nonumber \\ &+\frac{6^{-2}t_3t_5^2t_8-4t_1^6t_2t_5-2t_1^6t_3t_4+5t_1^4t_3t_8+4t_2t_5t_6+2t_3t_4t_6}{2(t_1^6-t_6)}\,\frac{\partial}{\partial
t_5} \nonumber
\\&+(-6t_2t_6)\,\frac{\partial}{\partial
t_6}+\frac{6^{-2}t_4^2-t_1^2}{2\times
6^{-2}}\,\frac{\partial}{\partial
t_7}+\frac{-3t_1^6t_2t_8+3t_1^5t_3t_8+3t_2t_6t_8}{t_1^6-t_6}\,\frac{\partial}{\partial
t_8}\,,
 \nonumber\\
\rvf&=t_1\frac{\partial}{\partial t_1}+2t_2\frac{\partial}{\partial
t_2}+3t_3\frac{\partial}{\partial t_3}+t_4\frac{\partial}{\partial
t_4}+2t_5\frac{\partial}{\partial
t_5}+6t_6\frac{\partial}{\partial t_6}+3t_8\frac{\partial}{\partial t_8}\,, \label{eq rvf4}\\
\cvf&=\frac{\partial}{\partial t_2}\,, \label{eq cvf4}
\end{align}}
where $t_8^2=36(t_1^6-t_6)$.
\end{itemize}
\end{exam}
Next, we state the main theorem of \cite{movnik} which, in a more general context, is also treated in \cite{younes1}.

\begin{theo}
 \label{main3}
There exist a unique vector field $\Ra=\Ra_n \in \mathfrak{X}(\T)$,
and unique regular functions $\Yuk_i\in \O_\T, \ 1\leq i\leq n-2,$
 such that:
\begin{equation}
\label{jimbryan} \nabla_{\Ra}
\begin{pmatrix}
\alpha_1\\
\alpha_2 \\
\alpha_3 \\
\vdots \\
\alpha_n \\
\alpha_{n+1} \\
\end{pmatrix}
= \underbrace {\begin{pmatrix}
0& 1 & 0&0&\cdots &0&0\\
0&0& \Yuk_1&0&\cdots   &0&0\\
0&0&0& \Yuk_2&\cdots   &0&0\\
\vdots&\vdots&\vdots&\vdots&\ddots   &\vdots&\vdots\\
0&0&0&0&\cdots   &\Yuk_{n-2}&0\\
0&0&0&0&\cdots   &0&-1\\
0&0&0&0&\cdots   &0&0\\
\end{pmatrix}}_\Yuk
\begin{pmatrix}
\alpha_1\\
\alpha_2 \\
\alpha_3 \\
\vdots \\
\alpha_n \\
\alpha_{n+1} \\
\end{pmatrix},
\end{equation}
 and $\Yuk\imc+\imc \Yuk^\tr=0$.
\end{theo}
Note that the vector field $\Ra=\Ra_n$ introduced in the above theorem is called modular vector field and, referring to the notation introduced in \eqref{eq A_H}, we have $\gm_\Ra=\Yuk$. In the
case where $n=3$, $\Ra_3$ is explicitly computed in \cite{ho22} and it is
verified that $\Yuk_1$ is the Yukawa coupling introduced in
\cite{can91}. $\Ra_1, \Ra_2, \Ra_4$ are explicitly described in
\cite{movnik} where the solutions of $\Ra_1$ and $\Ra_2$ are in
terms of (quasi-)modular forms, and when $n=4$ we found that
$\Yuk_1^2=\Yuk_2^2$ is the same as $4$-point function presented in
\cite[Table~1, $d=4$]{grmorpl}.

For any positive integer $n$, let $\LG:=\LG_n$ be the following algebraic
group
\begin{equation} \label{eq LG}
\LG:=\{\gG\in \textrm{GL}(n+1,\C)\, | \ \gG \rm{\ is  \ upper \
triangular \ and \ } \gG^\tr \Phi \gG=\Phi \},
\end{equation}
where $\Phi$ is the constant matrix given in \eqref{eq phi odd}. Let us fix the notation $\di2:=\frac{n+1}{2}$ if
$n$ is an odd positive integer, and $\di2:=\frac{n}{2}$ if $n$ is an
even positive integer.
$\LG$ acts on $\T$ from the right and it has $\di2$ copies of
$\mathbb{G}_m$ and $\dt-(\di2+1)$ copies of $\mathbb{G}_a$ as
subgroups where $\mathbb{G}_m$ and $\mathbb{G}_a$, respectively,
refer to the multiplicative group $(\C^\ast,\cdot)$ and the additive
group $(\C,+)$. We present an interpretation of
this group action in the chart $\t=(t_1,t_2,\ldots,t_\dt)$ of $\T$.
We observe that in the case of $n=1$ this action coincides with the
action of $\mathbb{G}_m$ and $\mathbb{G}_a$ given in \cite[\S
6.3]{ho14} for the family of elliptic curves arising from
Weierstrass form of elliptic curves, and in the case $n=3$ this
action has been studied in \cite{ho22} (see Section \ref{section Lie
group}).

We see that the Lie algebra of $\LG$ is
\begin{equation} \label{eq LA}
\LA=\{\gL\in{\rm Mat}(n+1,\C)\ | \ \gL \ {\rm is \ upper \
triangular \ and} \ \gL^\tr\Phi+\Phi\gL=0  \}.
\end{equation}
$\LA$ is a $\dt-1$ dimensional Lie algebra and we find its canonical
basis formed by $\gL_{\a \b}$'s, $1\leq \a \leq \di2, \ \a\leq
\b \leq 2\di2 +1-\a$, given in \eqref{eq g_ab1} and \eqref{eq
g_ab2}. In the following theorem we illustrate that to every element
of $\LA$ there corresponds a unique vector field on $\T$.

\begin{theo} \label{theo 1}
For any $\gL\in\LA$, there exists a unique vector field $\Ra_\gL\in
\mathfrak{X}(\T)$ such that:
\begin{equation}\label{29jul2017}
\gm_{\Ra_\gL}=\gL^\tr,
\end{equation}
i.e., $\nabla_{\Ra_\gL}\alpha =\gL^\tr\alpha$.
\end{theo}

The Lie algebra generated by  $\Ra_{\gL_{\a \b}}$'s, $1\leq \a \leq
\di2, \ \a\leq \b \leq 2\di2 +1-\a$, in $\mathfrak{X}(\T)$ with the Lie
bracket of the vector fields is isomorphic to $\LA$ with the Lie bracket of
the matrices. Hence, for simplicity of the notation, we use $\LA$ alternately either as a Lie subalgebra
of $\mathfrak{X}(\T)$ or as a Lie subalgebra of ${\rm Mat}(n+1,\C)$. Note that
$\Ra_{\gL_{\a\b}}$'s construct the canonical basis of $\LA\subseteq \mathfrak{X}(\T)$.

We define the \emph{AMSY-Lie
algebra} $\amsy$ to be the $\O_\T$-module generated by $\LA$
and the modular vector field $\Ra$ in $\mathfrak{X}(\T)$. To determine the Lie algebra structure of $\amsy$, it is enough to find
the Lie brackets $[\Ra,\Ra_{\gL_{\a\b}}]$ which is done in the following theorem. Before stating the theorem we fix some notations.
In what follows, $\delta_j^k$ denotes the Kronecker
delta, $\varrho(n)=1$ if $n$ is an odd integer and $\varrho(n)=0$
if $n$ is an even integer, $\Yuk_j$'s, $1\leq j\leq n-2$, are the functions given in Theorem \ref{main3}, and besides them we let
$\Yuk_0=-\Yuk_{n-1}=1$.

\begin{theo} \label{theo 2}
Considering the notations defined above, the followings hold:
 \begin{align}
&{[\Ra,\Ra_{\gL_{11}}]}=\Ra, \label{eq Rag11}\\
&{[\Ra,\Ra_{\gL_{22}}]}=-\Ra, \label{eq Rag22} \\
&{[\Ra,\Ra_{\gL_{\a\a}}]}=0, \ 3\leq \a \leq m\, , \label{eq Ragaa} \\
&[\Ra,{\Ra_{\gL_{\a\b}}}]=\Psi_1^{\a\b}(\Yuk)\,{\Ra_{\gL_{(\a+1)\b}}}+\Psi_2^{\a\b}(\Yuk)\,{\Ra_{\gL_{\a(\b-1)}}},\,
 1\leq \a \leq \di2, \ \a+1\leq \b \leq 2\di2+1-\a \, , \label{eq Ragab}
\end{align}
where
\begin{align}
&\Psi_1^{\a\b}(\Yuk):=(1+\varrho(n)\delta_{\a+\b}^{2\di2}-\delta_{\a+\b}^{2\di2+1})\,\Yuk_{\a-1},\\
&\Psi_2^{\a\b}(\Yuk):= (1-2\varrho(n)\delta_{\b}^{\di2+1})\,\Yuk_{n+1-\b}\,.
\end{align}
\end{theo}

If $n=1,2$, then we observe that $\amsy$ is isomorphic to
$\mathfrak{sl}_2(\C)$. In general, for $n\geq 3$ we
find a copy of $\mathfrak{sl}_2(\C)$ containing $\Ra$ as a Lie subalgebra of $\amsy$, which is stated in the following
theorem.

\begin{theo} \label{theo 3} Let the vector fields $\rvf$ and $\cvf$ be defined as follows:
\begin{description}
  \item[(i)] if $n=1$, then $\rvf:=-\Ra_{\gL_{11}}$ and $\cvf:=\Ra_{\gL_{12}}$,
  \item[(ii)] if $n=2$, then $\rvf:=-2\Ra_{\gL_{11}}$ and $\cvf:=2\Ra_{\gL_{12}}$,
  \item[(iii)] if $n\geq 3$, then $\rvf:=\Ra_{\gL_{22}}-\Ra_{\gL_{11}}$ and $\cvf:=\Ra_{\gL_{12}}$.
\end{description}
Then the Lie algebra generated by the vector fields
$\Ra, \rvf, \cvf$ in $\mathfrak{X}(\T)$ is isomorphic to
$\mathfrak{sl}_2(\C)$; indeed we get:
\[
[\Ra,\cvf]=\rvf \ , \ \ [\rvf,\Ra]=2\Ra \ , \ \ [\rvf,\cvf]=-2\cvf
\, .
\]
\end{theo}
In $\O_\T$, for $n=1,2,3,4$,  we attach to any $t_i$ the weight $\deg(t_i):=w_i$,
where $w_i$'s are defined by $\rvf=\sum_{i=1}^\dt
w_it_i\frac{\partial}{\partial t_i}$. If we let the vector fields $\Ra$ and $\cvf$ operate on $\O_{\T}$ as derivations, then $\cvf$ is a
derivation of degree $-2$, and the modular vector field $\Ra$ is a derivation of degree $2$ on the weighted space $\O_{\T}$ (see Examples \ref{exam n=1,4} and \ref{exam n=1,2}). Note that
a derivation $D$ on $\O_{\T}$ is of degree $r$ if for a given quasi-homogenous function $\varphi\in \O_{\T}$ we have $\deg(D\varphi)=\deg(\varphi)+r$. The author believes that these facts hold for any positive integer $n$.

\begin{rem} \label{rem cymf}
It is worth briefly mentioning a possible application of Theorems \ref{theo 2} and \ref{theo 3}.
On one hand, it is well known that the space of modular forms $M_\ast$ for $SL(2,\Z)$ is given by $\C[E_4,E_6]$, i.e., $M_\ast$ is generated by the Eisenstein series $E_4,E_6$. On the other hand, the triple $(E_2,E_4,E_6)$ of the Eisenstein series gives a solution of the Ramanujan's vector field \eqref{eq ram vf}. Hence by following the same process introduced by Don Zagier in \cite[\S5]{zag94} for the Eisenstein series $E_2,E_4,E_6$ one may define the Serre derivation and the Rankin-Cohen bracket on the space of \cy modular forms (see \cite{younes3}).
\end{rem}

The outline of this paper is as follows: Section \ref{section MSGMCD} contains a brief summary of
the relevant facts and the terminologies of \cite{movnik} which are necessary to have a self-contained manuscript. In Section \ref{section Lie group} we
discuss the action of the algebraic group $\LG$ on $\T$; indeed we observe that $\dim (\LG)=\dt-1$ and we
present its action in the chart $\t=(t_1,t_2,\ldots,t_\dt)$ of $\T$. Section \ref{section Lie algebra} is devoted to introduce $\LA$
and AMSY-Lie algebra $\amsy$. In this section we also prove Theorem \ref{theo 1} and Theorem \ref{theo 2}.
In Section \ref{subsection sl2} we construct a copy of $\mathfrak{sl}_2(\C)$ as a Lie subalgebra of $\amsy$ which contains the modular vector field $\Ra$ and we demonstrate Theorem \ref{theo 3}. The last section, Section \ref{subsection cy3}, establishes a brief summary of the AMSY-Lie algebra $\amsy$ attached to a non-rigid compact \cy threefold $X$ and shows that $\amsy$ contains $h^{21}$ copies of $\sl2$ where $h^{21}$ is the Hodge number of type $(2,1)$ of $X$.\\

{\bf Acknowledgment.} The author would like to appreciate Hossein Movasati
for valuable conversations that benefited him during the preparation of this paper. He also
thanks Julio Rebelo for commenting the importance of $\mathfrak{sl}_2(\C)$ in their works which
encouraged the author more to continue with the present research. During the preparation period of this manuscript the author was partially supported by
"Funda\c{c}\~ao Carlos Chagas Filho de Amparo \`{a} Pesquisa do Estado do Rio de Janeiro (FAPERJ)"
with process number E-26/010.001735/2016.
\section{Moduli spaces and GMCD}\label{section MSGMCD}
For the convenience of the reader and also to have a
self-contained text, in this section we recall some relevant facts
and terminologies discussed in \cite{movnik}, which construct the foundation of the present article.

In the Dwork family \eqref{12jan2016} it is more convenient to substitute the parameter $\psi$ with the parameter $z:=\psi^{-(n+2)}$ which is considered
as the standard parameter in the literatures. Let $W_z$ be the $n$-dimensional hypersurface in $\P^{n+1}$ given by:
\[
f_z(x_0,x_1,\ldots,x_{n+1}):=zx_0^{n+2}+x_1^{n+2}+x_2^{n+2}+\cdots+x_{n+1}^{n+2}-(n+2)x_0
x_1x_2\cdots x_{n+1}=0.
\]
$W_z$ represents a one parameter family of \cy $n$-folds. The finite group $G:=\{(\zeta_0,\zeta_1,\ldots,\zeta_{n+1})\mid \,\, \zeta_i^{n+2}=1,
\ \zeta_0\zeta_1\ldots \zeta_{n+1}=1 \}$,
acts canonically on $W_z$ as
$$
(\zeta_0,\zeta_1,\ldots,\zeta_{n+1}).(x_0,x_1,\ldots,x_{n+1})=(\zeta_0x_0,\zeta_1x_1,\ldots,\zeta_{n+1}x_{n+1}).
$$
The mirror variety $X=X_z,\ z\in \Pn 1\setminus\{0,1,\infty\}$, is obtained by the desingularization of the quotient
space $W_z/G $ (for more details see \cite[\S2]{movnik}).

By the \emph{moduli space of holomorphic $n$-forms} $\Ts$ we mean the
moduli of the pairs $(X,\alpha_1)$ where $X$ is an $n$-dimensional
mirror variety and $\alpha_1$ is a holomorphic $n$-form on $X$. We
know that the family of mirror varieties $X_z$ is a one parameter
family and the $n$-form $\alpha_1\,$ is unique, up to multiplication
by a constant, therefore $\dim( \Ts)=2$. Analogously to the construction of $X_z$, let $\X_{t_1,t_{n+2}}$,
$(t_1,t_{n+2})\in \C^2 \setminus \{(t_1^{n+2}-t_{n+2})t_{n+2}=0\}$,
be the mirror variety obtained by the desingularization
of the quotient space of the variety arising from the equation
\[
f_{t_1,t_{n+2}}(x_0,x_1,\ldots,x_{n+1}):=t_{n+2}x_0^{n+2}+x_1^{n+2}+x_2^{n+2}+\cdots+x_{n+1}^{n+2}-(n+2)t_1x_0
x_1x_2\cdots x_{n+1}=0.
\]
We fix two $n$-forms $\eta$ and $\omega_1$, respectively, in the
families $X_z$ and $\X_{t_1,t_{n+1}}$ such that in the
affine space $\{x_0=1\}$ are given as follows:
\begin{equation}
\eta:=\frac{dx_1\wedge dx_2\wedge \ldots \wedge dx_{n+1}}{df_z} \ ,
\ \ \  \omega_1:=\frac{dx_1\wedge dx_2\wedge \ldots \wedge
dx_{n+1}}{df_{t_1,t_{n+2}}}\ .
\end{equation}
Any element of $\Ts$ is in the form $(X_z,a\eta)$ where $a$ is a
non-zero constant. $(X_z,a\eta)$ can be identified by
$(\X_{t_1,t_{n+2}},\omega_1)$ as follows:
\begin{align}
&(X_z,a\eta)\mapsto
   (\X_{t_1,t_{n+2}},\omega_1)\, , \qquad
   (t_1,t_{n+2})=(a^{-1},za^{-(n+2)}) \, ,\\
   &(\X_{t_1,t_{n+2}},\omega_1)\mapsto
   (X_z,t_1^{-1}\eta) \, ,\qquad z=\frac{t_{n+2}}{t_1^{n+2}}\, .
\end{align}
Hence $(t_1,t_{n+2})$ construct a chart for $\Ts$; in the other word
$$\Ts=\spec ( \C[t_1,t_{n+2},\frac{1}{(t_1^{n+2}-t_{n+2})t_{n+2}}])\,,$$
and the morphism $\X\to\Ts$ is the universal family of the moduli of the pairs
$(X,\alpha_1)$. The multiplicative group $\Gm:=(\C^*,\cdot)$ acts on
$\Ts$ from the right by $ (X,\alpha_1)\bullet k=(X,k^{-1}\alpha_1)$ where $ k\in
\Gm,\ (X,\alpha_1)\in \Ts$. This action can be interpreted in the chart
$(t_1,t_{n+2})$ as:
\begin{equation}
\label{poloar} (t_1,t_{n+2})\bullet k=(kt_1,k^{n+2}t_{n+2}),\
(t_1,t_{n+2})\in \Ts,\ k\in \Gm \,,
\end{equation}
which follows from the isomorphism
\begin{align}
 &(\X_{kt_1,k^{n+2}t_{n+2}},\ k\omega_1)\cong (\X_{t_1,t_{n+2}}, \omega_1), \label{20m2014}\\
 &(x_1,x_2,\cdots ,x_{n+1})\mapsto
 (k^{-1}x_1,k^{-1}x_2,\cdots,k^{-1}x_{n+1})\, .\nonumber
\end{align}
Let
$
 \nabla:H_{\dR}^{n}(\X/\Ts)\to \Omega_\Ts^1\otimes_{\O_\Ts}H_{\dR}^{n}(\X/\Ts)
$ be the Gauss-Manin connection of the  two parameter family of
varieties $\X/\Ts$.  We define the $n$-forms $\omega_i,\,\
i=1,2,\ldots, n+1$, as follows:
\begin{equation}
\label{29oct11} \omega_i:= (\nabla_{\frac{\partial}{\partial
t_1}})^{i-1}(\omega_1),
\end{equation}
in which $\frac{\partial}{\partial t_1}$ is considered as a vector
field on the moduli space $\Ts$.
Then $\omega:=\{\omega_1,\omega_2,\ldots,\omega_{n+1}\} $ forms a
basis of $H^{n}_\dR(X)$ which is compatible with its Hodge
filtration, i.e.,
\begin{equation}
\label{eq. Gr. tr.}\omega_i\in F^{n+1-i}\setminus F^{n+2-i},
i=1,2,\ldots, n+1,
\end{equation}
where $F^i$ is the $i$-th piece of the Hodge filtration of
$H^n_\dR(X)$.
We can write the Gauss-Manin connection of $\X/\Ts$ as follow
\begin{equation} \label{eq Atilde}
\nabla\omega=\gma \omega\,, {\ \rm with \ } \omega={{\left( {\begin{array}{*{20}{c}}
  {{\omega _1}}&{{\omega _2}}& \ldots &{{\omega _{n + 1}}}
\end{array}} \right)}^{tr}}.
\end{equation}
If we denote by $\gma[i,j]$ the $(i,j)$-th entry of the Gauss-Manin
connection matrix $\gma$, then we obtain:
\begin{align}
&\gma[i,i] = -\frac{i}{(n+2)t_{n+2}}dt_{n+2}\, , \,\ 1 \leq i \leq n\, \label{16/1/2016-1}, \\
&\gma[i,i+1]= dt_1-\frac{t_1}{(n+2)t_{n+2}}dt_{n+2}\, , \,\ 1 \leq i \leq n\, , \label{16/1/2016-2}  \\
&\gma[n+1,j]=\frac{-S_2(n+2,j)t_1^j}{t_1^{n+2}-t_{n+2}}dt_1+\frac{S_2(n+2,j)t_1^{j+1}}{(n+2)t_{n+2}(t_1^{n+2}-t_{n+2})}dt_{n+2}\, , \,\ 1 \leq j \leq n\, , \label{eq stir1} \\
&\gma[n+1,n+1]=\frac{-S_2(n+2,n+1)t_1^{n+1}}{t_1^{n+2}-t_{n+2}}dt_1+
\frac{\frac{n(n+1)}{2}t_1^{n+2}+(n+1)t_{n+2}}{(n+2)t_{n+2}(t_1^{n+2}-t_{n+2})}dt_{n+2}\,, \label{eq stir2}
\end{align}
where $S_2(r,s)$ is the Stirling number of the second kind defined
by \begin{equation} \label{16jan2016} S_2\left( {r,s} \right){\rm{
}} = {\rm{ }}\frac{1}{{s!}}{\rm{ }}\sum\limits_{i = 0}^s {{{( -
1)}^i} \left( {\begin{array}{*{20}{c}}
s\\
i
\end{array}} \right)} {\left( {s - i} \right)^r}\, ,
\end{equation}
and the rest of the entries of $\gma$ are zero. As we mentioned in \cite{movnik}, equations \eqref{eq stir1} and \eqref{eq stir2} are checked for  $n=1,2,3,4$ and we believe that they are valid for arbitrary $n$, even though such explicit expressions do not interfere with our proofs. For any
$\xi_1,\xi_2\in H^n_\dR(X)$, in the context of the de Rham cohomology,
the \emph{intersection form} of $\xi_1$ and $\xi_2$, denoted by
$\langle \xi_1,\xi_2 \rangle$, is given as
$$
\langle \xi_1,\xi_2 \rangle:=\frac{1}{(2\pi i)^n}\int_{X}\xi_1\wedge
\xi_2\,,
$$
which is a non-degenerate $(-1)^n$-symmetric form. We obtain
\begin{align}
  & \langle \omega_i,\omega_j \rangle =0, {\rm \ if} \  i+j\leq n+1\, , \\
  &\langle \omega_1,\omega_{n+1} \rangle
  =(-(n+2))^n\frac{c_n}{t_1^{n+2}-t_{n+2}}, {\rm \ where} \ c_n \  {\rm is \
  a \  constant}\, ,\\
  &\langle \omega_j,\omega_{n+2-j} \rangle =(-1)^{j-1}\langle \omega_1,\omega_{n+1} \rangle, {\rm \ for}  \ j=1,2,\ldots,n+1\, ,
\end{align}
from which we can determine all the rest of
$\langle\omega_i,\omega_j\rangle$'s in a unique way. If we set
$\Omega=\Omega_n:=\left( \langle
\omega_i,\omega_j\rangle\right)_{1\leq i,j\leq n+1}, $ to be the
intersection form matrix in the basis $\omega$, then we have
\begin{equation}
\label{14/12/2015} d\Omega=\gma \Omega+ \Omega \gma^{\tr}.
\end{equation}
For any positive integer $n$ by the enhanced moduli space $\T=\T_n$ we mean the moduli of the pairs $(X,[\alpha_1,\cdots
,\alpha_n,\alpha_{n+1}])$ where $X$ is an $n$-dimensional mirror variety and $\{\alpha_1,\alpha_2,\ldots,\alpha_{n+1}\}$ constructs a basis of $H^n_\dR(X)$
satisfying the properties
$$
\alpha_i\in F^{n+1-i}\setminus F^{n+2-i},\ \ i=1,\cdots,n,n+1,
$$
and
\begin{equation}\label{11jan2016}
[\langle \alpha_i,\alpha_j\rangle]_{1\leq i,j \leq n+1}=\imc_n.
\end{equation}
Here, $\imc=\imc_n$ is the following constant $(n+1)\times(n+1)$
matrix:
\begin{equation}\label{eq phi odd}
\imc_n:=\left( {\begin{array}{*{20}c}
   {0_{\di2} } & {J_{\di2}}   \\
   { - J_{\di2}}  & {0_{\di2}}   \\
\end{array}} \right) \, {\rm if \, }n {\rm \, is \, odd,\, and} \,
\imc_n:=J_{n+1} \, {\rm if \, }n {\rm \, is \, even,}
\end{equation}
where by $0_{k}, k\in \mathbb{N},$ we mean a $k\times k$ block of
zeros, and $J_k$ is the following $k\times k$ block
\begin{equation}\small
J_{k }  := \left( {\begin{array}{*{20}c}
   0 & 0 &  \ldots  & 0 & 1  \\
   0 & 0 &  \ldots  & 1 & 0  \\
    \vdots  &  \vdots  &  {\mathinner{\mkern2mu\raise1pt\hbox{.}\mkern2mu
 \raise4pt\hbox{.}\mkern2mu\raise7pt\hbox{.}\mkern1mu}}  &  \vdots  &  \vdots   \\
   0 & 1 &  \ldots  & 0 & 0  \\
   1 & 0 &  \ldots  & 0 & 0  \\
\end{array}} \right).
\end{equation}
 We find that
\begin{equation}
\label{29dec2015} \dt=\dt_n:=\dim ( \T)=\left \{
\begin{array}{l}
\frac{(n+1)(n+3)}{4}+1,\,\, \quad  \textrm{\rm if \textit{n} is odd}
\\\\
\frac{n(n+2)}{4}+1,\,\,\,\,\quad\quad \textrm{\rm if \textit{n} is
even}
\end{array} \right. .
\end{equation}
Next, we are going to present a chart for the enhanced moduli space $\T$. In order to do
this, let $S=\left(
                           \begin{array}{c}
                             s_{ij} \\
                           \end{array}
                         \right)_{1\leq i,j\leq n+1}$ be a lower triangular matrix whose  entries are indeterminates $s_{ij},\ \ i\geq j$ and $s_{11}=1$.
We define
\[\underbrace {{{\left( {\begin{array}{*{20}{c}}
  {{\alpha _1}}&{{\alpha _2}}& \ldots &{{\alpha _{n + 1}}}
\end{array}} \right)}^{tr}}}_\alpha  = S\underbrace {\,\,{{\left( {\begin{array}{*{20}{c}}
  {{\omega _1}}&{{\omega _2}}& \ldots &{{\omega _{n + 1}}}
\end{array}} \right)}^{tr}}}_\omega \, ,\]
which implies that $\alpha$  forms a basis of $H^n_\dR(X)$
compatible with its Hodge filtration. We would like that
$(X,[\alpha_1,\alpha_2,\ldots,\alpha_{n+1}])$ be a member of $\T$,
hence it has to satisfy $\left(
                           \begin{array}{c}
                             \langle\alpha_i,\alpha_j\rangle \\
                           \end{array}
                         \right)_{1\leq i,j\leq n+1}=\Phi$, from what we get the following equation
\begin{equation}\label{eq sost}
S\Omega S^\tr=\Phi.
\end{equation}
Using this equation we can express $\frac{(n+2)(n+1)}{2}-\dt-2$
numbers of the parameters $s_{ij}$'s in terms of the other $\dt-2$
parameters that we fix them as \emph{independent parameters}. For
simplicity we write the first class of the parameters as $\check
t_1,\check t_2,\cdots, \check t_{d_0}$ and the second class, which are the independent parameters, as $t_2,
t_3,\ldots,t_{n+1},t_{n+3},\ldots,t_\dt$. We put all these
parameters inside $S$ according to the following rule which we write
it only for $n=1,2,3,4,5$:
\[\tiny \left( {\begin{array}{*{20}{c}}
  1&0 \\
  {{t_2}}&{{{\check t}_1}}
\end{array}} \right) , \left( {\begin{array}{*{20}{c}}
1&0&0\\
{{t_2}}&{{\check t}_2}&0\\
{{\check t}_4}&{{\check t}_3}&{{\check t}_1}
\end{array}} \right),\left( {\begin{array}{*{20}{c}}
1&0&0&0\\
{{t_2}}&{{t_3}}&0&0\\
{{t_4}}&{{t_6}}&{{\check t}_2}&0\\
{{t_7}}&{{\check t}_4}&{{\check t}_3}&{{\check t}_1}
\end{array}} \right), \left( {\begin{array}{*{20}{c}}
1&0&0&0&0\\
{{t_2}}&{{t_3}}&0&0&0\\
{{t_4}}&{{t_5}}&{{\check t}_3}&0&0\\
{{t_7}}&{{\check t}_7}&{{\check t}_5}&{{\check t}_2}&0\\
{{\check t}_9}&{{\check t}_8}&{{\check t}_6}&{{\check t}_4}&{{\check
t}_1}
\end{array}} \right), \left( {\begin{array}{*{20}{c}}
  1&0&0&0&0&0 \\
  {{t_2}}&{{t_3}}&0&0&0&0 \\
  {{t_4}}&{{t_5}}&{{t_6}}&0&0&0 \\
  {{t_8}}&{{t_9}}&{{t_{10}}}&{{{\check t}_3}}&0&0 \\
  {{t_{11}}}&{{t_{12}}}&{{{\check t}_7}}&{{{\check t}_5}}&{{{\check t}_2}}&0 \\
  {{t_{13}}}&{{{\check t}_9}}&{{{\check t}_8}}&{{{\check t}_6}}&{{{\check t}_4}}&{{{\check t}_1}}
\end{array}} \right)
.\]
Note that we have already used  $t_1,t_{n+2}$ as coordinates
system of $\Ts$. In particular we find:
\begin{equation}
\label{29/12/2015} s_{(n+2-i)(n+2-i)}=
\frac{(-1)^{n+i+1}}{c_n(n+2)^n}\frac{t_1^{n+2}-t_{n+2}}{s_{ii}}, \
1\leq i \leq m \,.
\end{equation}
Hence $\t:=(t_1,t_2,\ldots,t_\dt)$ forms a chart for the enhanced moduli space
$\T$, and in fact
\begin{align}
 \label{thanksdeligne}
 \T&=\spec(\C[t_1,t_2,\ldots, t_{\dt},\frac{1}{t_{n+2}(t_{n+2}-t_1^{n+2})\check t }])\,,\\
 \O_\T&=\C[t_1,t_2,\ldots, t_{\dt},\frac{1}{t_{n+2}(t_{n+2}-t_1^{n+2})\check t }]\,,
\end{align}
Here, $\check t$ is the product of $\di2-1$ independent variables which are located on the main diagonal of $S$.
From now on, we alternately use either $s_{ij}$'s, or $t_i$'s
and ${\check t}_j$'s to refer the entries of $S$. If we denote by
$\gm$  the Gauss-Manin connection matrix of the family $\X/\T$
written in the basis $\alpha$, i.e., $\nabla \alpha=\gm \alpha$,
then we calculate $\gm$ as follow:
\begin{equation}\label{eq gm 3/27/2016}
\gm=\left (dS+S\cdot \gma\right )\, S^{-1}\,.
\end{equation}
In the following remarks we recall some results deduced from the
proof of Theorem \ref{main3} in \cite[\S 7]{movnik}.

\begin{rem}
We obtain the functions $\Yuk_i$'s given in \eqref{jimbryan} as
follows: if $n$ is odd then
\begin{align}
&\Yuk_{i}=-\Yuk_{n-(i+1)}=\frac{s_{22}\,s_{(i+1)(i+1)}}{s_{(i+2)(i+2)}},\
\ i=1,2,\ldots, \frac{n-3}{2}\, , \label{eq yukiodd}\\
&\Yuk_{\frac{n-1}{2}}=(-1)^{\frac{3n+3}{2}}c_n(n+2)^n\frac{s_{22}\,s_{\frac{n+1}{2}\frac{n+1}{2}}^2}{t_1^{n+2}-t_{n+2}}\,
, \label{eq yukmodd}
\end{align}
and if $n$ is even then
\begin{align}
&\Yuk_{i}=-\Yuk_{n-(i+1)}=\frac{s_{22}\,s_{(i+1)(i+1)}}{s_{(i+2)(i+2)}},\
\ i=1,2,\ldots, \frac{n-2}{2}\, . \label{eq yukieven}
\end{align}
\end{rem}
\begin{rem}\label{rem uniqueness}
Let $\H\in \mathfrak{X}(\T)$. If $\ \nabla_\H\alpha=0$ for any
$(X,[\alpha_1,\alpha_2,\ldots,\alpha_{n+1}])\in \T$, then $\H=0$.
\end{rem}
\section{Algebraic group}\label{section Lie group}
For any positive integer $n$, let $\LG:=\LG_n$ be the algebraic
group given in \eqref{eq LG}.
$\LG$ acts on $\T$ from the right as follow:
\[
(X,[\alpha_1,\alpha_2,\ldots,\alpha_{n+1}])\bullet
\gG=(X,\alpha^\tr\gG),
\]
where $\alpha:=\left( {\begin{array}{*{20}{c}} {{\alpha
_1}}&{{\alpha _2}}& \ldots &{{\alpha _{n + 1}}}
\end{array}} \right)^\tr$, $\gG\in\LG$,  and  in the right hand side of the above equality
$\alpha^\tr\gG$ refers to the matrix product. It is of interest to
interpret this action in the chart $\t$ presented for $\T$. To this
end, we first announce some properties of $\LG$. Before, note
that we fixed $m:=\frac{n+1}{2}$ if $n$ is an odd integer, and
$m:=\frac{n}{2}$ if $n$ is an even integer.

The equation $\gG^\tr \Phi \gG=\Phi \ $ in \eqref{eq LG} guaranties
that $\dim ( \LG)=\dt-1$. For $i=1,2,\ldots,m$, by defining
\[\LG_{i}:=\left\{ {\gG=\left(\gG_{kl}\right)\in
\LG\,\left| \begin{array}{l}  \ \gG {\rm \ is \ a \ diagonal \
matrix \ with\ }
\gG_{ii}=\gG_{(n+2-i)(n+2-i)}^{-1}\in \C^\ast,\\
 \ {\rm and} \  \gG_{kk}=1 \ {\rm for \ all} \ k\neq
i,n+2-i .
\end{array} \right.} \right\}\]
we find out that $\LG_i\simeq \mathbb{G}_m$; hence $\LG$ has $m$ copies of $\mathbb{G}_m$ as
multiplicative subgroups. To describe additive subgroups of $\LG$,
we consider the following two cases:
\begin{description}
  \item[First case:] If $n$ is odd, then for $1\leq i \leq \frac{n+1}{2}$ and $i+1\leq j \leq n+2-i$ define $\LG_{ij}$ as follow:
\[\LG_{ij}:=\left\{ {\gG=\left(\gG_{kl}\right)\in
\LG\,\left| \begin{array}{l}
\blacktriangleright \ \gG_{kk}=1,\ {\rm for} \  k=1,2,\ldots,n+1,\\
\blacktriangleright \ {\rm if} \ j\leq m, \ {\rm then} \ \gG_{ij}=-\gG_{(n+2-j)(n+2-i)}\in \C,\\
\blacktriangleright \ {\rm if} \ j\geq m+1, \ {\rm then} \ \gG_{ij}=\gG_{(n+2-j)(n+2-i)}\in \C, \\
\blacktriangleright \ {\rm and \ the \ rest \ of \ the \ entries \
are \ zero}.
\end{array} \right.} \right\}.\]
One can easily check that $\LG_{ij}\simeq \mathbb{G}_a$, i.e., $\LG_{ij}$ is an additive
subgroup of $\LG$. Thus we find $\dt-(m+1)$($=\frac{(n+1)^2}{4}$) copies of $\mathbb{G}_a$ as
additive subgroups of $\LG$.
\item[Second case:] If $n$ is even, then for $1\leq i \leq \frac{n}{2}$ and $i+1\leq j \leq n+1-i$ we consider $\LG_{ij}$ as follow:
\[\LG_{ij}:=\left\{ {\gG=\left(\gG_{kl}\right)\in
\LG\,\left| \begin{array}{l}
\blacktriangleright \ \gG_{kk}=1,\ {\rm for} \  k=1,2,\ldots,n+1,\\
\blacktriangleright \ \gG_{ij}=-\gG_{(n+2-j)(n+2-i)}\in \C,\\
\blacktriangleright \ {\rm if} \ j= \frac{n+2}{2}, \ {\rm then} \
\gG_{i(n+2-i)}=-\frac{1}{2}\gG_{(n+2-j)(n+2-i)}^2 , \\
\blacktriangleright \ {\rm and \ the \ rest \ of \ the \ entries \
are \ zero}.
\end{array} \right.} \right\}.\]
Again it is not difficult to show that $\LG_{ij}\simeq \mathbb{G}_a$. Therefore in this case we have
$\dt-(m+1)$($=\frac{n^2}{4}$) copies of additive subgroups of $\LG$ as well.
\end{description}
We give the new order $\LG_{m+1}\ldots,\LG_{\dt-1}$ to additive
subgroups $\LG_{ij}$'s, in such a way that in this order
$\LG_{i_1j_1}$ appears before $\LG_{i_2j_2}$ provided that $i_1<i_2$, and in
the case that $i_1=i_2$ then $\LG_{i_1j_1}$ appears before $\LG_{i_1j_2}$ if $j_1<j_2$. Any $\gG_i\in\LG_i, \
i=1,2,\ldots,\dt-1$, can be presented by a unique complex number, that by abuse of notation, we denote it again by $\gG_i$. For any $\gG\in\LG$ there are
unique elements $\gG_i\in\LG_i,i=1,2,\ldots,\dt-1$, such that
$\gG=\gG_1\gG_2\ldots\gG_{\dt-1}$. Therefore, we can represent any
$\gG\in\LG$ by a $(\dt-1)$-tuple $(\gG_1,\gG_2,\ldots,\gG_{\dt-1})$,
where $\gG_i\in \C^\ast$ provided $i=1,2,\ldots,m$, and $\gG_i\in \C$ provided
$i=m+1,m+2,\ldots,\dt-1$; in the other words
\[\LG \simeq \underbrace {\mathbb{G}_m \times \mathbb{G}_m \times  \ldots  \times \mathbb{G}_m}_{m {\rm - times}}
\times \underbrace {\mathbb{G}_a \times \mathbb{G}_a \times  \ldots
\times \mathbb{G}_a}_{\dt - (m + 1){\rm - times}}.\] Hence we can
summarize the above facts in the following lemma.
\begin{lemm}\label{lemm LG}
$\LG$ is a $(\dt-1)$-dimensional Lie group which contains $m$ copies of
$\mathbb{G}_m$ as multiplicative subgroups and $\dt-(m+1)$ copies
of $\mathbb{G}_a$ as additive subgroups. Moreover, for any
$\gG\in\LG$, there exist unique elements
$\gG_i\in \LG_i,\ i=1,2,\ldots,\dt-1$, such that
$\gG=\gG_1\gG_2\ldots\gG_{\dt-1}$.
\end{lemm}

The following example helps to have a better imagination of elements
of $\LG_i$'s.

\begin{exam}
If n=3, then any $\gG\in \LG$ can be written  as
$\gG=\gG_1\gG_2\gG_3\gG_4\gG_5\gG_6$ where $\gG_i\in \LG_i$ are
given as follows:{\tiny
\[\begin{array}{l}
{\gG_1} = \left( {\begin{array}{*{20}{c}}
{\gG_1^{ - 1}}&0&0&0\\
0&1&0&0\\
0&0&1&0\\
0&0&0&{{\gG_1}}
\end{array}} \right)\,,\,\,\,{\gG_2} = \left( {\begin{array}{*{20}{c}}
1&0&0&0\\
0&{\gG_2^{ - 1}}&0&0\\
0&0&{{\gG_2}}&0\\
0&0&0&1
\end{array}} \right)\,,\,\,\,{\gG_3} = \left( {\begin{array}{*{20}{c}}
1&{ - {\gG_3}}&0&0\\
0&1&0&0\\
0&0&1&{{\gG_3}}\\
0&0&0&1
\end{array}} \right)\\
{\gG_4} = \left( {\begin{array}{*{20}{c}}
1&0&{{\gG_4}}&0\\
0&1&0&{{\gG_4}}\\
0&0&1&0\\
0&0&0&1
\end{array}} \right)\,,\,\,\,{\gG_5} = \left( {\begin{array}{*{20}{c}}
1&0&0&{{\gG_5}}\\
0&1&0&0\\
0&0&1&0\\
0&0&0&1
\end{array}} \right)\,,\,\,\,\,{\gG_6} = \left( {\begin{array}{*{20}{c}}
1&0&0&0\\
0&1&{{\gG_6}}&0\\
0&0&1&0\\
0&0&0&1
\end{array}} \right)\,.
\end{array}\]
} If n=4, then any $\gG\in \LG$ can be written  as
$\gG=\gG_1\gG_2\gG_3\gG_4\gG_5\gG_6$, in which $\gG_i\in \LG_i$ are
given bellow:{\tiny
\[\begin{array}{l}
{\gG_1} = \left( {\begin{array}{*{20}{c}}
{\gG_1^{ - 1}}&0&0&0&0\\
0&1&0&0&0\\
0&0&1&0&0\\
0&0&0&1&0\\
0&0&0&0&{{\gG_1}}
\end{array}} \right)\,,\,\,\,{\gG_2} = \left( {\begin{array}{*{20}{c}}
1&0&0&0&0\\
0&{\gG_2^{ - 1}}&0&0&0\\
0&0&1&0&0\\
0&0&0&{{\gG_2}}&0\\
0&0&0&0&1
\end{array}} \right)\,,\,\,\,{\gG_3} = \left( {\begin{array}{*{20}{c}}
1&{ - {\gG_3}}&0&0&0\\
0&1&0&0&0\\
0&0&1&0&0\\
0&0&0&1&{{\gG_3}}\\
0&0&0&0&1
\end{array}} \right)\\
{\gG_4} = \left( {\begin{array}{*{20}{c}}
1&0&{ - {\gG_4}}&0&{ - \frac{1}{2}\gG_4^2}\\
0&1&0&0&0\\
0&0&1&0&{{\gG_4}}\\
0&0&0&1&0\\
0&0&0&0&1
\end{array}} \right)\,,\,\,\,{\gG_5} = \left( {\begin{array}{*{20}{c}}
1&0&0&{ - {\gG_5}}&0\\
0&1&0&0&{{\gG_5}}\\
0&0&1&0&0\\
0&0&0&1&0\\
0&0&0&0&1
\end{array}} \right)\,,\,\,\,\,{\gG_6} = \left( {\begin{array}{*{20}{c}}
1&0&0&0&0\\
0&1&{ - {\gG_6}}&{ - \frac{1}{2}\gG_6^2}&0\\
0&0&1&{{\gG_6}}&0\\
0&0&0&1&0\\
0&0&0&0&1
\end{array}} \right).
\end{array}\]
}
\end{exam}

From \eqref{20m2014} we get that for any $k\in\C^\ast$
\begin{equation}\label{eq iso ms}
  (\X_{t_1,t_{n+2}},S\omega) \cong (\X_{kt_1,k^{n+2}t_{n+2}},S{\tiny \left( {\begin{array}{*{20}{c}}
k&0&\ldots&0&0\\
0&{{k^2}}&\ldots&0&0\\
\vdots&\vdots& \ddots &\vdots&\vdots\\
0&0& \ldots &k^n&0\\
0&0&\ldots&0&{{k^{n + 1}}}
\end{array}} \right)}\omega),
\end{equation}
in which $\omega=\left( {\begin{array}{*{20}{c}} {{\omega
_1}}&{{\omega _2}}& \ldots &{{\omega _{n + 1}}}
\end{array}} \right)^\tr$ and $S$ is the basis change matrix
$\alpha=S\omega$.

Let $(\X_{t_1,t_{n+2}},[\alpha_1,\alpha_2,\ldots,\alpha_{n+1}])$ be
a presentation of $\t=(t_1,t_2,\ldots,t_\dt)\in\T$. Then for any
$\gG=(\gG_1,\gG_2,\ldots,\gG_{n+1})\in\LG$ we have:
\begin{align}
(\X_{t_1,t_{n+2}},[\alpha_1,\alpha_2,\ldots,\alpha_{n+1}])\bullet\gG
&=(\X_{t_1,t_{n+2}},(S\omega)^\tr \gG) \nonumber \\
&\cong (\X_{\gG_1t_1,\gG_1^{n+2}t_{n+2}},\left(S{\tiny \left( {\begin{array}{*{20}{c}}
\gG_1&0&\ldots&0&0\\
0&{{\gG_1^2}}&\ldots&0&0\\
\vdots&\vdots& \ddots &\vdots&\vdots\\
0&0& \ldots &\gG_1^n&0\\
0&0&\ldots&0&{{\gG_1^{n + 1}}}
\end{array}} \right)}\omega\right)^\tr\gG)\ ,\label{eq last pair}
\end{align}
in which the latter isomorphism is concluded from \eqref{eq iso ms}. The last pair of the above equations, namely \eqref{eq last pair}, gives the action $(t_1,t_2,\ldots,t_\dt)\bullet \gG$. In particular, if we denote by
$(t_1,t_2,\ldots,t_{\dt})\bullet\gG=(t_1\bullet\gG,t_2\bullet\gG,\ldots,t_{\dt}\bullet\gG)$, then we always have $t_1\bullet\gG=t_1\gG_1$ and $t_{n+2}\bullet\gG=t_{n+2}\gG_1^{n+2}$.
\begin{exam}
For $n=1,2,3,4$ we state the action of
$\gG=(\gG_1,\gG_2,\ldots,\gG_{\dt-1})\in\LG$ on
$\t=(t_1,t_2,\ldots,t_{\dt}) \in\T$ as follows:
\begin{description}
  \item[n=1:]
  \[
    (t_1,t_2,t_3)\bullet \gG=(t_1\gG_1,t_2\gG_1^2+\gG_2,t_3\gG_1^3).
  \]
  \item[n=2:]
  \[
    (t_1,t_2,t_4)\bullet \gG=(t_1\gG_1,t_2\gG_1-\gG_2,t_4\gG_1^4).
  \]
  \item[n=3:]
\begin{align}
   &t_1\bullet\gG=t_1\gG_1, \qquad\qquad\qquad t_2\bullet\gG=(t_2\gG_1-\gG_2\gG_3)\gG_2^{-1}, \nonumber \\
   &t_3\bullet\gG=t_3\gG_1^2\gG_2^{-1}, \ \qquad\qquad t_4\bullet\gG=(t_2\gG_1\gG_6+t_4\gG_1\gG_2^2-\gG_2\gG_3\gG_6+\gG_2\gG_4)\gG_2^{-1}, \nonumber \\
   &t_5\bullet\gG=t_5\gG_1^5, \qquad\qquad\qquad t_6\bullet\gG=(t_3\gG_1^2\gG_6+t_6\gG_1^2\gG_2^2)\gG_2^{-1}, \nonumber \\
   &t_7\bullet\gG=(t_2\gG_1\gG_4+t_4\gG_1\gG_2^2\gG_3+t_7\gG_1^2\gG_2-\gG_2\gG_3\gG_4+\gG_2\gG_5)\gG_2^{-1}\,. \nonumber
  \end{align}

  \item[n=4:]
   \begin{align}
   &t_1\bullet\gG=t_1\gG_1, \qquad\qquad\qquad\qquad\quad\ \ t_2\bullet\gG=(t_2\gG_1-\gG_2\gG_3)\gG_2^{-1}, \nonumber \\
   &t_3\bullet\gG=t_3\gG_1^2\gG_2^{-1}, \nonumber \qquad\qquad\qquad\qquad t_4\bullet\gG=(-t_2\gG_1\gG_6+t_4\gG_1\gG_2+\gG_2\gG_3\gG_6-\gG_2\gG_4)\gG_2^{-1}, \nonumber \\
   &t_5\bullet\gG=(-t_3\gG_1^2\gG_6+t_5\gG_1^2\gG_2)\gG_2^{-1}, \qquad t_6\bullet\gG=t_6\gG_1^6, \nonumber \\
   &t_7\bullet\gG=\frac{1}{2}(-t_2\gG_1\gG_6^2+2t_4\gG_1\gG_2\gG_6+2t_7\gG_1\gG_2^2+\gG_2\gG_3\gG_6^2-2\gG_2\gG_4\gG_6-2\gG_2
\gG_5)\gG_2^{-1}, \nonumber \\
   &t_8\bullet\gG=t_8\gG_1^3\,. \nonumber
  \end{align}
\end{description}
\end{exam}

\section{Lie algebra}\label{section Lie algebra}

We stated the Lie algebra $\LA$ of $\LG$ in \eqref{eq LA}.
If $n$ is an odd integer, then for
$1\leq \a \leq \di2$ and  $\a\leq \b \leq n+2-\a$ let
\begin{equation} \label{eq g_ab1}
\gL_{\a\b}:={\left( {{g_{kl}}} \right)_{(n + 1) \times (n + 1)}} , \
{\rm such \ that \ } \left\{ {\begin{array}{*{20}{c}}
{{\rm if \ } \b\leq m, {\rm \ then \ }{g_{\a\b}} = 1,\,\,{g_{(n + 2 - \b)(n + 2 - \a)}} =  - 1,}\\
{{\rm if \ } \b\geq m+1, {\rm \ then \ }{g_{\a\b}} ={g_{(n + 2 - \b)(n + 2 - \a)}} =  1,}\ \\
{{\rm and \ the \ rest \ of \ the \ entries \ are \ zero.
}\qquad\qquad \ \, }
\end{array}} \right.
\end{equation}
If $n$ is an even integer, then for $1\leq \a \leq \di2$ and $\a\leq \b
\leq n+1-\a$ set
\begin{equation} \label{eq g_ab2}
\gL_{\a\b}:={\left( {{g_{kl}}} \right)_{(n + 1) \times (n + 1)}} , \
{\rm such \ that \ } \left\{ {\begin{array}{*{20}{c}}
{{g_{\a\b}} = 1,\,\,{g_{(n + 2 - \b)(n + 2 - \a)}} =  - 1,}\qquad \ \ \\
{{\rm and \ the \ rest \ of \ the \ entries \ are \ zero. } \, }
\end{array}} \right.
\end{equation}
One can easily check that the set of $\gL_{\a\b}$'s forms a canonical
basis of $\LA$.

\begin{exam}
If n=3, then \eqref{eq g_ab1} yields:{\tiny
\[\begin{array}{l}
{\gL_{11}} = \left( {\begin{array}{*{20}{c}}
1&0&0&0\\
0&0&0&0\\
0&0&0&0\\
0&0&0&-1
\end{array}} \right)\,,\,\,\,{\gL_{22}} = \left( {\begin{array}{*{20}{c}}
0&0&0&0\\
0&1&0&0\\
0&0&-1&0\\
0&0&0&0
\end{array}} \right)\,,\,\,\,{\gL_{12}} = \left( {\begin{array}{*{20}{c}}
0&1&0&0\\
0&0&0&0\\
0&0&0&-1\\
0&0&0&0
\end{array}} \right)\\
{\gL_{13}} = \left( {\begin{array}{*{20}{c}}
0&0&1&0\\
0&0&0&1\\
0&0&0&0\\
0&0&0&0
\end{array}} \right)\,,\,\,\,{\gL_{14}} = \left( {\begin{array}{*{20}{c}}
0&0&0&1\\
0&0&0&0\\
0&0&0&0\\
0&0&0&0
\end{array}} \right)\,,\,\,\,\,{\gL_{23}} = \left( {\begin{array}{*{20}{c}}
0&0&0&0\\
0&0&1&0\\
0&0&0&0\\
0&0&0&0
\end{array}} \right)\,.
\end{array}\]
} If n=4, then by \eqref{eq g_ab2} we have:{\tiny
\[\begin{array}{l}
{\gL_{11}} = \left( {\begin{array}{*{20}{c}}
1&0&0&0&0\\
0&0&0&0&0\\
0&0&0&0&0\\
0&0&0&0&0\\
0&0&0&0&-1
\end{array}} \right)\,,\,\,\,{\gL_{22}} = \left( {\begin{array}{*{20}{c}}
0&0&0&0&0\\
0&1&0&0&0\\
0&0&0&0&0\\
0&0&0&-1&0\\
0&0&0&0&0
\end{array}} \right)\,,\,\,\,{\gL_{12}} = \left( {\begin{array}{*{20}{c}}
0&1&0&0&0\\
0&0&0&0&0\\
0&0&0&0&0\\
0&0&0&0&-1\\
0&0&0&0&0
\end{array}} \right)\\
{\gL_{13}} = \left( {\begin{array}{*{20}{c}}
0&0&1&0&0\\
0&0&0&0&0\\
0&0&0&0&-1\\
0&0&0&0&0\\
0&0&0&0&0
\end{array}} \right)\,,\,\,\,{\gL_{14}} = \left( {\begin{array}{*{20}{c}}
0&0&0&1&0\\
0&0&0&0&-1\\
0&0&0&0&0\\
0&0&0&0&0\\
0&0&0&0&0
\end{array}} \right)\,,\,\,\,\,{\gL_{23}} = \left( {\begin{array}{*{20}{c}}
0&0&0&0&0\\
0&0&1&0&0\\
0&0&0&-1&0\\
0&0&0&0&0\\
0&0&0&0&0
\end{array}} \right).
\end{array}\]
}
\end{exam}

\subsection{Proof of Theorem \ref{theo 1}}
This proof is somehow analogous to the proof of Theorem \ref{main3}
given in \cite[\S 7]{movnik}. We first construct another moduli
space $\tilde \T$ which contains $\T$. In order to do this, assume
that all the entries $s_{ij},\ j\leq i, (i,j)\not=(1,1)$, of $S$ are
independent parameters. We denote by $\tilde\T$ and $\tilde\alpha$
the corresponding family of varieties and a basis of differential
forms. Indeed, $\tilde \T$ is formed the same as $\T$ by removing the
condition \eqref{11jan2016}.

Let $\gL\in \LA$ be arbitrary. In account of \eqref{eq gm
3/27/2016}, the existence of a vector field
\[\Ra_\gL:={\dot
t}_1\frac{\partial}{\partial t_1}+{\dot
t}_{n+2}\frac{\partial}{\partial t_{n+2}}+\sum_{i=2,j=1}^{n+1,i}
{\dot s}_{ij}\frac{\partial}{\partial s_{ij}}\,,
\]
in $\tilde\T$ with the desired property \eqref{29jul2017} is
equivalent to the existence of a solution of the equation
\begin{equation}
\label{eq dots2} \dot{S}=\gm_{\Ra_\gL}\cdot S-S\cdot
\gma(\Ra_\gL).
\end{equation}
Note that here $\dot x:=dx(\Ra_\gL)$ is the derivation of the
function $x$ along the vector field  $\Ra_\gL$ in $\tilde\T$. In \eqref{eq dots2} the
equality corresponding to the $(i,j)$-th entry, $j\leq i,\ \
(i,j)\not=(1,1)$, serves as the definition of $\dot{s}_{ij}$. The
equality corresponding to $(1,1)$-th and $(1,2)$-th entries of \eqref{eq dots2} yield,
respectively, $\dot{t}_1$ and $\dot{t}_{n+2}$. All the rest of the entries of \eqref{eq dots2} are
trivial equalities $0=0$. Therefore, we conclude the statement of Theorem
\ref{main3} for the moduli space $\tilde\T$.

To prove the statement for the moduli space $\T$, consider the map
\begin{equation}
 \tilde\T\to \Mat_{(n+1)\times(n+1)}(\C),\ \ (t_1,t_{n+2}, S)\mapsto S \Omega
 S^{\tr}\,.
\end{equation}
It follows that $\T$ is the fiber of this map over the point $\Phi$.
It is enough to prove that the vector field $\Ra_\gL$ is tangent to the fiber of
the above map over $\Phi$. This follows from
\begin{eqnarray*}
\overbrace{(S\, \Omega\, S^{\tr})}^{.}  &=&
\dot S\,\Omega \, S^{\tr}+S\,\dot\Omega\, S^{\tr}+ S\,\Omega\, \dot S^{\tr}\\
&=&
(\gm_{\Ra_\gL} S-S\,{\gma}(\Ra_\gL) )\,\Omega\, S^{\tr}+S \,({\gma}(\Ra_\gL)\, \Omega+\Omega\,\gma^{\tr}(\Ra_\gL))\, S^{\tr}+ S\Omega ( S^{\tr}\gm_{\Ra_\gL}^{\tr}-\gma^{\tr}(\Ra_\gL) S^{\tr})\\
&=&
\gm_{\Ra_\gL} \Phi+ \Phi\,\gm_{\Ra_\gL}^{\tr} \\
&=&
\gL^\tr \Phi+ \Phi\,\gL \\
&=& 0\,,
\end{eqnarray*}
where $\dot x:=dx(\Ra_\gL)$ is the derivation of the function $x$
along the vector field $\Ra_\gL$ in $\T$.
 Note that in the above equalities we are using \eqref{14/12/2015} and
the fact that $\gL$ belongs to $\LA$.

The uniqueness of $\Ra_\gL$ follows from Remark \ref{rem
uniqueness}.

\subsection{$\LA$ as a Lie subalgebra of $\mathfrak{X}(\T)$}

We know that the Gauss-Manin connection $\nabla$ is a flat (or
integrable) connection, i.e., for any $\H_1,\H_2\in
\mathfrak{X}(\T)$ we have
$\nabla_{[\H_1,\H_2]}=\nabla_{\H_1}\nabla_{\H_2}-\nabla_{\H_2}\nabla_{\H_1}$
where $[\H_1,\H_2]$ refers to the Lie bracket of the vector fields. If $\nabla_{\H_1}\alpha=\gm_{\H_1}\alpha$ and $\nabla_{\H_2}\alpha=\gm_{\H_2}\alpha$, then one obtains
\[
\nabla_{\H_1}\nabla_{\H_2}\alpha=\gm_{\H_2}\gm_{\H_1}\alpha+\H_1(\gm_{\H_2}).
\]
Hence, by letting
$\nabla_{[{\H_1},{\H_2}]}\alpha=\gm_{[\H_1,\H_2]}\alpha$, we
get:
\begin{equation} \label{eq liebracamsy}
\gm_{[\H_1,\H_2]}=[\gm_{\H_2},\gm_{\H_1}]+
\H_1(\gm_{\H_2})-\H_2(\gm_{\H_1}),
\end{equation}
where on the right hand side by $[\gm_{\H_2},\gm_{\H_1}]$ we mean
the Lie bracket of the matrices. Since elements of $\LA$ are constant
matrices, for any  $\gL_1,\gL_2\in \LA$ we obtain:
\begin{equation}\label{eq liebrackLA}
\gm_{[\Ra_{\gL_{1}},\Ra_{\gL_{2}}]}=[\gm_{\Ra_{\gL_{2}}},\gm_{\Ra{\gL_{1}}}]\,.
\end{equation}
Consider the map
\[
\varphi: \, \LA \to \mathfrak{X}(\T)\,, \ \ \ \varphi(\gL)=\Ra_\gL
\,.
\]
For any $\gL_1,\gL_2\in \LA$, by using of \eqref{29jul2017}
 we have:
\begin{align}
&\gm_{\varphi([\gL_1,\gL_2])}=[\gL_1,\gL_2]^\tr=[\gm_{\Ra_{\gL_2}},\gm_{\Ra_{\gL_1}}]\,.
\label{eq liebrackLie(G)}
\end{align}
Therefore, in account of \eqref{eq liebrackLA}, \eqref{eq
liebrackLie(G)} and Remark \ref{rem uniqueness} we get
$\varphi([\gL_1,\gL_2]=[\Ra_{\gL_1},\Ra_{\gL_2}]$ which implies
that the Lie algebra generated by $\Ra_{\gL_{\a\b}}$'s, $ 1\leq
\a\leq m, \, \a\leq \b\leq 2m+1-\a$ in $\mathfrak{X}(\T)$ is
isomorphic to $\LA$. By abuse of notation, we use the same notation $\LA$ for the
Lie algebra generated by $\Ra_{\gL_{\a\b}}$'s, and by employing
\eqref{eq liebrackLA} and Remark \ref{rem uniqueness} we can
determine the Lie bracket of $\LA\subset \mathfrak{X}(\T)$
completely.

\subsection{AMSY-Lie algebra}

As we mentioned in Section \ref{section introduction}, we define
AMSY-Lie algebra $\amsy$ as $\O_{\T}$-module
generated by $\LA$ and the modular vector field $\Ra$, which is a Lie subalgebra of $\mathfrak{X}(\T)$.
Determining the Lie algebra structure of $\amsy$ is not
as easy as of $\LA$,  since not all matrices corresponded to the elements of $\amsy$ are constant.
In order to do this, it is enough to determine
$[\Ra,\Ra_{\gL_{\a\b}}]$ for $1\leq \a\leq m, \, \a\leq \b\leq
2m+1-\a$. Before
proving Theorem \ref{theo 2}, note that if we denote by
$\H=\sum_{i=1}^{\dt}\dot{t}_{i} \frac{\partial}{\partial t_i}\in
\mathfrak{X}(\T)$, then we find $\dot{t}_{i}$'s from
\begin{equation} \label{eq Sdot}
\dot{S}=\gm_\H S-S\, \gma(\H)\,,
\end{equation}
where $\dot x:=dx(\H)$ is the derivation of the function $x$ along
the vector field $\H$ in $\T$. We saw in Section \ref{section
MSGMCD} that any $t_i,\,1\leq i\leq \dt,\, i\neq 1,n+2$, corresponds
to only one $s_{jk}$, $1\leq j,k\leq n+1$.

\subsection{Proof of Theorem \ref{theo 2}}

Equation \eqref{eq liebracamsy} yields:
\begin{equation}\label{eq liebracamsy2}
\gm_{[\Ra,\Ra_{\gL_{\a\b}}]}=[\gm_{\Ra_{\gL_{\a\b}}},\Yuk]-\Ra_{\gL_{\a\b}}(\Yuk).
\end{equation}
Hence, to determine $[\Ra,\Ra_{\gL_{\a\b}}]$, we need to compute
$\Ra_{\gL_{\a\b}}(\Yuk)$. If we apply \eqref{eq Sdot} to
$\Ra_{\gL_{11}}$, then the equalities corresponding to $(1,1)$-th
and $(1,2)$-th entries, respectively, lead to
$$
\dot{t}_1=-t_1,\ \ \dot{t}_{n+2}=-(n+2)t_2,
$$
from what we obtain the diagonal matrix $\gma(\Ra_{\gL_{11}})=\diag(1,2,\ldots,n+1)$. This implies
\[
\Ra_{\gL_{11}}=\sum_{i=1}^{\dt}c_{i\gL_{11}} t_i
\frac{\partial}{\partial t_i}\ ,
\]
in which $c_{i\gL_{11}}$'s are constants given as follows:
\begin{equation}\label{eq coef Rg1}
c_{i\gL_{11}}= \left\{ {\begin{array}{*{20}{c}}
{{-1}, \ \ {\rm if}} \ \ i=1 \ ,\qquad\qquad\qquad\qquad\quad \ \ \ \ \, \qquad\qquad\qquad\qquad\qquad\quad \\
{{-(n+2), \ \ {\rm if}} \ \ i=n+2 \ ,\qquad\qquad\quad \ \ \ \, \qquad\qquad\qquad\qquad\qquad\quad}\\
{-2 , \ \ {\rm if} \ \ i=\dt \ \ {\rm and} \ \  n \ \ {\rm is \ \
odd}
\ ,\qquad \ \ \ \ \qquad\qquad\qquad\qquad\qquad\quad}\\
{-1 , \ \ {\rm if} \ \ i=\dt \ \ {\rm and} \ \  n \ \ {\rm is \ \
even}
\ ,\qquad \ \ \ \, \qquad\qquad\qquad\qquad\qquad\quad}\\
{-k, \ \ {\rm if} \ \ t_i=s_{jk}\ \ {\rm and} \ \  i\neq 1, n+2,
\dt \ {\rm for \ some}\ 2\leq j,k\leq n+1 . \ \ }
\end{array}} \right.
\end{equation}
For $2\leq \a \leq m$, we employ \eqref{eq Sdot} for
$\Ra_{\gL_{\a\a}}$ that yields $\dot t _1=0$ and $\dot t_{n+2}=0$.
Hence we find $\gma(\Ra_{\gL_{\a\a}})=0$ which provides
\[
\Ra_{\gL_{\a\a}}=\sum_{i=1}^{\dt}c_{i\gL_{\a\a}} t_i
\frac{\partial}{\partial t_i}\ ,
\]
where $c_{i\gL_{\a\a}}$'s are following constants:
\begin{equation}\label{eq coef Rgj}
c_{i\gL_{\a\a}}= \left\{ {\begin{array}{*{20}{c}}
{{1}, \ \ {\rm if} \ \ t_i=s_{\a k},\ {\rm for \ \ some \ \ }1\leq k \leq n+1 \ , \qquad }  \\
{{-1}, \ \ {\rm if} \ \ t_i=s_{(n+2-\a)k},\ {\rm for \ \ some \ \ }1\leq k \leq n+1 \ ,}\\
{0, \ \ {\rm otherwise}\ .
\qquad\qquad\qquad\qquad\qquad\qquad\qquad }
\end{array}} \right.
\end{equation}
Analogously, for $\Ra_{\gL_{\a\b}}$, $1\leq \a\leq m, \, \a+1\leq
\b\leq 2m+1-\a$, we obtain $\gma(\Ra_{\gL_{\a\b}})=0$, thus  $\Ra_{\gL_{\a\b}}$ follows from the equation
\begin{equation}\label{eq Rgij}
\dot{S}=\gm_{\Ra_{\gL_{\a\b}}} S \ .
\end{equation}
The equations \eqref{eq coef Rg1} and \eqref{eq coef Rgj} yield:
\begin{align}
&\Ra_{\gL_{11}}=-(t_1\frac{\partial}{\partial
t_1}+(n+2)t_{n+2}\frac{\partial}{\partial
t_{n+2}})-\sum_{k=2}^{m}ks_{kk}\frac{\partial}{\partial
s_{kk}}+{ \widehat{\Ra}_{\gL_{11}}},\\
&\Ra_{\gL_{\a\a}}=s_{\a\a}\frac{\partial}{\partial s_{\a\a}}+{ \widehat{\Ra}_{\gL_{\a\a}}}\, , \ \ \a=2,3,\ldots, m\, ,
\end{align}
where in the part $\widehat{\Ra}_{\gL_{\a\a}}, \, \a=1,2,3,\ldots, m$, do not appear the terms including
$t_1\frac{\partial}{\partial t_1}$, $t_{n+2}\frac{\partial}{\partial
t_{n+2}}$ and $s_{jj}\frac{\partial}{\partial s_{jj}} , \
j=2,3,\ldots,m$. Hence in account of equations \eqref{eq yukiodd},
\eqref{eq yukmodd} and \eqref{eq yukieven} we get:
\begin{align}
&\Ra_{\gL_{11}}(\Yuk_i)=-\Yuk_i\ , \ \  1\leq i \leq m-1  ,\\
&\Ra_{\gL_{22}}(\Yuk_1)=2\Yuk_1, \ \Ra_{\gL_{22}}(\Yuk_i)=\Yuk_i, \
\
2\leq i \leq m-1, \\
&\Ra_{\gL_{\a\a}}(\Yuk_{\a-2})=-\Yuk_{\a-2}, \
\Ra_{\gL_{\a\a}}(\Yuk_{\a-1})=\Yuk_{\a-1}, \ 3\leq \a \leq m-1,\\
&\Ra_{\gL_{\a\a}}(\Yuk_{i})=0, \ 3\leq \a \leq m-1, \ 1\leq i \leq
m-1 \ {\rm and} \ i\neq \a-2,\a-1\,.
\end{align}
If $n$ is odd, then
\begin{equation}
\Ra_{\gL_{mm}}(\Yuk_{m-2})=-\Yuk_{m-2}, \
\Ra_{\gL_{mm}}(\Yuk_{m-1})=2\Yuk_{m-1}, \
\Ra_{\gL_{mm}}(\Yuk_{i})=0, \ 1\leq i \leq m-3,
\end{equation}
and if $n$ is even, then
\begin{equation}
\Ra_{\gL_{mm}}(\Yuk_{m-2})=-\Yuk_{m-2}, \
\Ra_{\gL_{mm}}(\Yuk_{m-1})=\Yuk_{m-1}, \ \Ra_{\gL_{mm}}(\Yuk_{i})=0,
 \ 1\leq i \leq m-3.
\end{equation}
These relations together with \eqref{eq liebracamsy2} and Remark
\ref{rem uniqueness} prove \eqref{eq Rag11}, \eqref{eq Rag22} and
\eqref{eq Ragaa}. We conclude from \eqref{eq Rgij}
that the terms including
$\ast\frac{\partial}{\partial s_{kk}}, \ 2\leq k \leq n+1,$ do not appear in the expression of $\Ra_{\gL_{\a\b}}$, $1\leq \a\leq m, \,
\a+1\leq \b\leq 2m+1-\a$. Hence $\Ra_{\gL_{\a\b}}(\Yuk)=0$ which combining with \eqref{eq
liebracamsy2} and Remark \ref{rem uniqueness} finishes the proof of
\eqref{eq Ragab}.

\section{$\mathfrak{sl}_2(\C)$ Lie algebra and the weights} \label{subsection sl2}

In the following example we observe that for $n=1,2$ the Lie algebra $\amsy$ is isomorphic to $\sl2$.
\begin{exam}\label{exam n=1,2}
If $n=1$, then we have:
\begin{align}
&\tiny \Yuk = \left( {\begin{array}{*{20}{c}}
  0&1 \\
  0&0
\end{array}\,} \right)\,\,\,\,,\,\,\,\,\,\,\,\,{\gL_{11}^\tr} = \left( {\begin{array}{*{20}{c}}
  1&0 \\
  0&{ - 1}
\end{array}\,} \right)\,\,\,\,,\,\,\,\,\,\,\,\,{\gL_{12}^\tr} = \left( {\begin{array}{*{20}{c}}
  0&0 \\
  1&0
\end{array}\,} \right)\,.
\end{align}
By letting $\rvf:=-\Ra_{\gL_{11}}$ and $\cvf:=\Ra_{\gL_{12}}$, we get $\Ra$, $\rvf$ and $\cvf$ respectively as given by \eqref{eq mvf1}, \eqref{eq rvf1} and \eqref{eq cvf1}. For $n=2$ we obtain:
\begin{align}
&\tiny \Yuk = \left( {\begin{array}{*{20}{c}}
  0&1&0 \\
  0&0&{ - 1} \\
  0&0&0
\end{array}\,} \right)\,\,\,\,,\,\,\,\,\,\,\,\,{\gL_{11}^\tr} = \left( {\begin{array}{*{20}{c}}
  1&0&0 \\
  0&0&0 \\
  0&0&{ - 1}
\end{array}\,} \right)\,\,\,\,,\,\,\,\,\,\,\,\,{\gL_{12}^\tr} = \left( {\begin{array}{*{20}{c}}
  0&0&0 \\
  1&0&{ 0} \\
  0&-1&0
\end{array}} \right)\,,\\
&\Ra=(t_3-t_1t_2)\frac{\partial}{\partial
t_1}+(2t_1^2-\frac{1}{2}t_2^2)\frac{\partial}{\partial
t_2}+(-2t_2t_3+8t_1^3)\frac{\partial}{\partial t_3}+(-4t_2t_4)\frac{\partial}{\partial t_4},\\
&\Ra_{\gL_{11}}=-t_1\frac{\partial}{\partial
t_1}-t_2\frac{\partial}{\partial
t_2}-2t_3\frac{\partial}{\partial t_3}-4t_4\frac{\partial}{\partial t_4},\\
&\Ra_{\gL_{12}}=\frac{\partial}{\partial t_2}\,,
\end{align}
where the polynomial equation $t_3^2=4(t_1^4-t_4)$ holds among
$t_i$'s. In this case we define $\rvf:=-2\Ra_{\gL_{11}}$ and
$\cvf:=2\Ra_{\gL_{12}}$. Hence for $n=1,2$ in account of the correspondences $\Ra\mapsto \mmat,\ \cvf\mapsto \cmat, \ \rvf\mapsto \rmat$ we observe that $\amsy$ is isomorphic to
$\mathfrak{sl}_2(\C)$.
\end{exam}
{\bf Proof of Theorem \ref{theo 3}.}
The proofs of {\bf (i)} and {\bf (ii)} are given in Example \ref{exam n=1,2}. For $n\geq 3$,  we define
$\rvf:=\Ra_{\gL_{22}}-\Ra_{\gL_{11}}$ and $\cvf:=\Ra_{\gL_{12}}$.
The equation \eqref{eq Rgij} yields:
\begin{equation}\label{eq cvf 2019}
\cvf=\left\{ {\begin{array}{*{20}{c}} {\frac{\partial }{{\partial
{t_2}}}\,,\,\,\,\,{\rm if}\,\,\,n\,\,{\rm is\,\,even},\,\,\,\,\,\,\,\,\,\,\,\,\,\,\,\,\ \ \ }\\
\\ {\frac{\partial }{{\partial {t_2}}}\, - {s_{n1}}\frac{\partial
}{{\partial {t_\dt}}},\,\,\,\,{\rm if}\,\,\,n\,\,{ \rm
is\,\,odd}\,.}
\end{array}} \right.
\end{equation}
If we consider $\Ra\mapsto \mmat,\ \cvf\mapsto \cmat, \ \rvf\mapsto \rmat$,  then it is an immediate result of Theorem \ref{theo 2} that
the Lie algebra generated by $\Ra,\, \rvf, \, \cvf$ in
$\mathfrak{X}(\T)$ is isomorphic to $\mathfrak{sl}_2(\C)$. \hfill\(\square\)\\

In the following example we state $\Ra,\, \rvf, \, \cvf$ for  $n=3,4$.

\begin{exam}
If $n=3$, then we get:{\small
\begin{align}
\Ra&=(t_3-t_1t_2)\frac{\partial}{\partial
t_1}+\frac{t_3^3t_4-5^4t_2^2(t_1^5-t_5)}{5^4(t_1^5-t_5)}\,\frac{\partial}{\partial
t_2}+\frac{t_3^3t_6-3\times5^4t_2t_3(t_1^5-t_5)}{5^4(t_1^5-t_5)}\,\frac{\partial}{\partial
t_3} \nonumber \\&+(-t_2t_4-t_7)\frac{\partial}{\partial
t_4}+(-5t_2t_5)\frac{\partial}{\partial
t_5}+(-t_2t_6-2t_3t_4+5^5t_1^3)\frac{\partial}{\partial
t_6}+(-5^4t_1t_3-t_2t_7)\frac{\partial}{\partial t_7}\,,
 \nonumber\\
\rvf&=t_1\frac{\partial}{\partial t_1}+2t_2\frac{\partial}{\partial
t_2}+3t_3\frac{\partial}{\partial t_3}+5t_5\frac{\partial}{\partial
t_5}+t_6\frac{\partial}{\partial t_6}+2t_7\frac{\partial}{\partial t_7}\,, \nonumber\\
\cvf&=\frac{\partial}{\partial t_2}-t_4\frac{\partial}{\partial
t_7}\,. \nonumber
\end{align}}
For $n=4$ the vector fields $\Ra$, $\rvf$ and $\cvf$ are respectively given by \eqref{eq mvf4}, \eqref{eq rvf4} and \eqref{eq cvf4}.
\end{exam}

\begin{rem}
For $n=2$ and $n=4$, to avoid the
second root of $\check t_2$ and $\check t_3$, in the coordinates chart $\t$ we are considering
extra variables $t_3:=\check t_2$ and $t_8:=\check t_3$ respectively. Because of this, in these cases, we are stating
a polynomial equation following from \eqref{29/12/2015}. An analog argument is applied to the enhanced moduli spaces
corresponding to all positive even integers $n$.
\end{rem}

For any $t_i$ in $\O_{\T}$ let $\deg(t_i)=w_i$ where $w_i$'s are defined by $\rvf=\sum_{i=1}^\dt w_i t_i\frac{\partial}{\partial t_i}$.
By letting $\Ra=\sum_{i=1}^\dt \dot t_i\frac{\partial}{\partial t_i}$, in account of the computations given for
$n=1,2,3,4$, we observe that $\dot t_i$'s are rational functions satisfying $\deg(\dot t_i)=w_i+2$.
If we consider modular vector field $\Ra$ as a derivation on the weighted space $\O_{\T}$,
then $\Ra$ increases the degree of any quasi-homogeneous regular function of $\O_{\T}$ by two.
For $n=1,2$ we found a
solution of $\Ra$ (see \cite{movnik})  in terms of (quasi-)modular forms. It is well known that
the usual derivation of the (quasi-)modular forms increases their weight by
two. This observation is another evidence which convince us to relate the modular vector field
$\Ra$ to the classical theory of (quasi-)modular forms.
The author believe that for any positive integer $n$ the modular vector field $\Ra$ is a derivation
of degree $2$, i.e., $\Ra$  increases the degree
of any quasi-homogeneous regular function of $\O_{\T}$ by two.

Equations \eqref{eq coef Rg1} and \eqref{eq coef Rgj}
yield that for any positive even integer $n$ the term $\ast\frac{\partial}{\partial s_{n1}}$ does not
appear in the expression of $\rvf$, hence the weight of $t_k$
corresponding to $s_{n1}$ is zero.  This fact combining with \eqref{eq cvf 2019} implies that $\cvf$ is a
vector field of degree zero and $\rvf$ is the radial vector field in the
weighted space $\O_\T$. This observation coincides, in some sense, with the
definition of Halphen type and semi-simple vector fields given in \cite{gui07,guireb}.

All these facts about weights, may motivate one to pursue the studying of the
dynamics of the modular vector field $\Ra$ in the weighted space $\O_\T$ with the assigned weights $\deg(t_i)=w_i$.

\begin{rem} \label{rem [fR,E]}
Let $f\in \O_\T$ be a quasi-homogenous polynomial of degree $k$. Then $\rvf(f)=k\,f$ which implies $[\rvf,f\Ra]=(k+2)f\Ra$. In particular if $f=t_1^{n+2}-t_{n+2}$, then we get
\begin{equation}\label{eq}
[f\Ra,\cvf]=f\rvf \ , \ \ [\rvf,f\Ra]=(n+4)f\Ra \ \, .
\end{equation}
Therefore, instead of working with $\Ra$, one possibly use $f\Ra$ which is a polynomial vector field and is more convenient to work with.
\end{rem}

\subsection{Another copy of $\mathfrak{sl}_2(\C)$} \label{subsection another sl2}
For $n=3$, if we set
{\small\begin{align}\label{a}
\H&:=(t_3-t_1t_2)\frac{\partial}{\partial
t_1}+(-t_2^2)\,\frac{\partial}{\partial
t_2}+(-3t_2t_3)\,\frac{\partial}{\partial
t_3}+(-t_2t_4-t_7)\frac{\partial}{\partial
t_4}\nonumber \\&+(-5t_2t_5)\frac{\partial}{\partial
t_5}+(-t_2t_6-2t_3t_4+5^5t_1^3)\frac{\partial}{\partial
t_6}+(-5^4t_1t_3-t_2t_7)\frac{\partial}{\partial t_7}\,,\nonumber
\end{align}}
then we have
\begin{align}
\gm_\H&=\left( {\begin{array}{*{20}{c}}
0&1&0&0\\
\frac{t_2t_3^2t_6-t_3^3t_4}{5^4(t_1^5-t_5)}&\frac{-t_3^2t_6}{5^4(t_1^5-t_5)}&\frac{t_3^3}{5^4(t_1^5-t_5)}&0\\
\frac{t_2t_3t_6^2-t_3^2t_4t_6}{5^4(t_1^5-t_5)}&\frac{-t_3t_6^2}{5^4(t_1^5-t_5)}&\frac{t_3^2t_6}{5^4(t_1^5-t_5)}&-1\\
\frac{-t_2^2t_3t_6^2+2t_2t_3^2t_4t_6-t_3^3t_4^2}{5^4(t_1^5-t_5)}&\frac{t_2t_3t_6^2-t_3^2t_4t_6}{5^4(t_1^5-t_5)}&\frac{-t_2t_3^2t_6+t_3^3t_4}{5^4(t_1^5-t_5)}&0
\end{array}} \right)  \\
&=\Yuk-\frac{\, t_6}{\, t_3}\,\Yuk_1\gm_{\Ra_{\gL_{11}}}+\frac{t_2t_6-t_3t_4}{t_3}\Yuk_1 \gm_{\Ra_{\gL_{12}}}
+\frac{t_2t_6^2-t_3t_4t_6}{t_3^2}\Yuk_1 \gm_{\Ra_{\gL_{13}}}\nonumber \\
&+\frac{-t_2^2t_6^2+2t_2t_3t_4t_6-t_3^2t_4^2}{t_3^2}\Yuk_1 \gm_{\Ra_{\gL_{14}}}-\frac{\, t_6^2}{\, t_3^2}\,\Yuk_1\gm_{\Ra_{\gL_{23}}} \nonumber,
\end{align}
from what we get $\H\in \amsy$.
 For $n=4$ let
{\small
\begin{align}
\H&:=(t_3-t_1t_2)\frac{\partial}{\partial
t_1}+(-t_2^2)\,\frac{\partial}{\partial
t_2}+(-3t_2t_3)\,\frac{\partial}{\partial
t_3} +(-t_2t_4)\,\frac{\partial}{\partial
t_4}\nonumber
\\&+(-2t_2t_5-t_3t_4)\,\frac{\partial}{\partial
t_5} +(-6t_2t_6)\,\frac{\partial}{\partial
t_6}+\frac{6^{-2}t_4^2-t_1^2}{2\times
6^{-2}}\,\frac{\partial}{\partial
t_7}+(-3t_2t_8)\,\frac{\partial}{\partial
t_8}\,.
 \nonumber
\end{align}}
Then we obtain the $(3,3)$-th entry of $\gm_\H$ as $\frac{-3t_1^5t_3}{t_1^6-t_6}$, which implies $\H\notin \amsy$. In the both cases $n=3,4$ it follows that the Lie algebra generated by  $\H,\, \rvf, \, \cvf$ is isomorphic to  $\mathfrak{sl}_2(\C)$.
Here $\H$ is a polynomial vector field following from $\Ra=\sum_{i=1}^\dt \dot t_i\frac{\partial}{\partial t_i}$ by discarding the non-polynomial terms of $\dot t_i,\, i=1,2,\ldots,\dt$. Note that, the same as $\Ra$, the vector field $\H$ is a degree $2$ derivation on
the weighted space $\O_{\T}$, and Remark \ref{rem [fR,E]} holds for $\H$ as well.

\section{$\mathfrak{sl}_2(\C)$ attached to \cy threefolds} \label{subsection cy3}

In this section we suppose that $X$ is a non-rigid compact \cy threefold on $\C$, and $\hn:=h^{21}$ is the Hodge number of type $(2,1)$ of $X$.
Here, we first give a brief summary of AMSY-Lie algebra $\amsy$ attached to the \cy threefold $X$ which is discussed in \cite{alimov}. Then we observe that there are $\hn$ copies of $\sl2$ in $\amsy$.
Note that $\dim( H_\dR^3(X))=2\hn+2$ and its Hodge filtration is as follow:
$$
0=F^4\subset F^3\subset F^2\subset F^1\subset F^0=H_\dR^3(X)\,.
$$
The enhanced moduli space $\T$
is the moduli of the pairs $(X,[\alpha_1, \alpha_2,\ldots
,\alpha_{2\hn+2}])$ where $X$ is as above and $\{\alpha_i\}_{i=1}^{2\hn+2}$ is a basis of $H_\dR^3(X)$ with the properties:
\begin{align}
\alpha_1\in F^{3},\ \alpha_1,\alpha_2,\ldots, \alpha_{\hn+1}\in F^2,\ &\alpha_{1},\alpha_{2},\ldots, \alpha_{2\hn+1}\in F^1,\ \alpha_{1},\alpha_{2},\ldots, \alpha_{2\hn+2}\in F^0\,,\nonumber\\
&[\langle \alpha_i,\alpha_j\rangle]=\imc\, ,\nonumber
\end{align}
in which $\imc$ is the following constant matrix:
\begin{equation}
\label{31aug10}
\Phi:=
\begin{pmatrix}
0&0&0&-1\\
0& 0&\mathbbm{1}_{\hn \times \hn}&0\\
0& -\mathbbm{1}_{\hn\times \hn}&0&0\\
 1&0&0&0
\end{pmatrix}.
\end{equation}
Here, it is used $(2\hn+2)\times (2\hn+2)$ block matrices according to the decomposition $2\hn+2=1+\hn+\hn+1$ and
$\mathbbm{1}_{\hn\times \hn}$ denotes the $\hn\times \hn$ identity matrix. The algebraic group
\begin{equation}\nonumber
\LG:=\left \{\gG\in \textrm{GL}(2\hn+2,\C)\mid \gG \text{ is block upper triangular and  }  \gG^\tr\imc\gG=\imc\ \ \right \}
\end{equation}
acts from the right on $\T$ and its Lie algebra is as follow
\begin{equation}\nonumber
\LA=
\left \{\gL\in {\rm Mat}(2\hn+2,\C)\mid \gL \text{ is block upper triangular and  }  \gL^\tr\imc+\imc\gL=0\ \ \right \}.
\end{equation}
By block triangular we mean triangular with respect to the partition
$2\hn+2=1+\hn+\hn+1$. One finds that
$$
\dim(\LG)=\frac{3\hn^2+5\hn+4}{2},\ \ \ \dim(\T)=\hn+\dim(\LG).
$$

In \cite{alimov} it is proved that there are unique modular vector fields $\Ra_k, \ k=1,2,\ldots,\hn$, on $\T$ and unique
regular functions $\Yukk_{ijk}\in \O_\T,\ \ i,j,k=1,2,\ldots,\hn$, which are symmetric in $i,j,k$, such that
\begin{equation}\nonumber
\gm_{\Ra_k}=\left(
\begin{array}{*{4}{c}}
0 & \delta^j_k& 0 & 0 \\
0 & 0 & \Yukk_{kij}& 0 \\
0 & 0 & 0 & \delta_k^i \\
0 & 0 & 0 & 0
\end{array} \right). \ \ \ \ \ \ \
\end{equation}
Also, for any $\gL\in \LA$ there is  a unique vector field $\Ra_{\gL}$  in $\T$
such that
\begin{equation}\nonumber
\gm_{\Ra_\gL}=\gL^{\tr}.
\end{equation}
Here the AMSY-Lie algebra $\amsy$ is the $\O_\T$-module generated by the
vector fields
\begin{equation}
\label{29apr2014}
\Ra_i,\ \   \Ra_{\gL},\ \textrm{where}\ i=1,2,\ldots,\hn,\ \textrm{and} \ \gL\in \LA .
\end{equation}
Followings form the canonical basis of $\LA$:
\begin{equation} \nonumber
\label{gofLie}
\tgtwo_{ab}:=\left(
\begin{array}{cccc}
 0 & 0 & 0 & 0 \\
 0 & 0 & 0 & 0 \\
 0& \frac{1}{2} (\delta^{i}_{a} \delta^{j}_{b} + \delta^{i}_{b} \delta^{j}_{a} ) & 0 & 0 \\
 0 & 0 & 0 & 0 \\
\end{array}
\right)^\tr ,\,
\tgone_{a}=
\left(
\begin{array}{cccc}
 0 & 0 & 0 & 0 \\
 0 & 0 & 0 & 0 \\
 -\delta^{i}_{a} & 0 & 0 & 0 \\
 0 & \delta_{a}^{j} & 0 & 0 \\
\end{array}
\right)^\tr,\,
\tgzero:=\left(
\begin{array}{cccc}
 0 & 0 & 0 & 0 \\
 0 & 0 & 0 & 0 \\
 0 & 0 & 0 & 0 \\
 -1 & 0 & 0 & 0 \\
\end{array}
\right)^\tr,\,
\end{equation}
$$
\kgone^{a}:= \left(
\begin{array}{cccc}
 0 & 0 & 0 & 0 \\
 \delta^{a}_{i} & 0 & 0 & 0 \\
 0 & 0 & 0 & 0 \\
 0 & 0 & \delta^{a}_{j} & 0 \\
\end{array}
\right)^\tr,\,  \ \
\ggtwo^{a}_{b} :=
 \left(
\begin{array}{cccc}
 0 & 0 & 0 & 0 \\
 0 & -\delta_{i}^{a} \delta^{j}_{b}   & 0 & 0 \\
 0 & 0 & \delta_{b}^{i} \delta_{j}^{a} & 0 \\
 0 & 0 & 0& 0 \\
\end{array}\right)^\tr,\,
\ggzero_0:=
\left(
\begin{array}{cccc}
 -1 & 0 & 0 & 0 \\
 0 & 0 & 0 & 0 \\
 0 & 0 & 0 & 0 \\
 0 & 0 & 0& 1 \\
\end{array}\right)^\tr.
$$
The Lie algebra structure of $\amsy$ is given by the following table.
{\tiny
\begin{equation} \nonumber
\label{liebrackettable}
\begin{array}{|c|c|c|c|c|c|c|c|}
\hline
&\Ra_{\ggzero_0} & \Ra_{\ggtwo_{c}^{d}}&\Ra_{\tgtwo_{cd}} & \Ra_{\tgone_{c}} &  \Ra_{\tgzero} &\Ra_{\kgone^{c}}&\Ra_{c} \\

\hline
 \Ra_{\ggzero_0} &0&0&0&-\Ra_{\tgone_{c}}&-2\Ra_{\tgzero}&-\Ra_{\kgone^{c}} &\Ra_{c}\\

\hline
\Ra_{\ggtwo_{b}^{a}}&0&0&-\delta_{c}^{a} \Ra_{\tgtwo_{bd}}- \delta_{d}^{a} \Ra_{\tgtwo_{bc}}&- \delta^{a}_{c} \Ra_{\tgone_{b}}&0&\delta^{c}_{b} \Ra_{\kgone^{a}} &-\delta^{a}_{c} \Ra_{b}\\

\hline
\Ra_{\tgtwo_{ab}} &0&\delta_{a}^{d} \Ra_{\tgtwo_{bc}}+\delta_{b}^{d} \Ra_{\tgtwo_{ac}}&0&0&0&\frac{1}{2} (\delta_{a}^{c} \Ra_{\tgone_{b}} + \delta_{b}^{c} \Ra_{\tgone_{a}} )&-\frac{1}{2} ( \Yukk_{cbd} \Ra_{\ggtwo_{a}^{d}}  +\Yukk_{acd} \Ra_{\ggtwo_{b}^{d}}    )\\

\hline
\Ra_{\tgone_{a}} &\Ra_{\tgone_{a}} &\delta^{d}_{a} \Ra_{\tgone_{c}}&0&0&0&2 \delta^{c}_{a}  \Ra_{\tgzero}&2 \Ra_{\tgtwo_{ac}} -\Yukk_{acd} \Ra_{\kgone^{d}} \\

\hline
\Ra_{\tgzero} &2 \Ra_{\tgzero} &0&0&0&0&0&\Ra_{\tgone_{c}}\\

\hline
\Ra_{\kgone^{a}} &\Ra_{\kgone^{a}}&-\delta^{a}_{c} \Ra_{\kgone^{d}}&-\frac{1}{2} (\delta_{c}^{a} \Ra_{\tgone_{d}} + \delta_{d}^{a} \Ra_{\tgone_{c}} )& -2\delta^{a}_{c} \Ra_{\tgzero}&0&0&-\delta_{c}^{a} \Ra_{\ggzero_0} +\Ra_{\ggtwo_{c}^{a}}\\

\hline
\Ra_{a}&-\Ra_{a}&\delta_{a}^{d} \Ra_{c}&
\frac{1}{2} ( \Yukk_{ade} \Ra_{\ggtwo_{c}^{e}}  +\Yukk_{ace} \Ra_{\ggtwo_{d}^{e}}    )&
-2 \Ra_{\tgtwo_{ac}} +\Yukk_{ace} \Ra_{\kgone^{e}} &
-\Ra_{\tgone_{a}}&
\delta_{a}^{c} \Ra_{\ggzero_0} -\Ra_{\ggtwo^{c}_{a}}&
0\\
\hline
\end{array}
\end{equation}
}
For $k=1,2,\ldots,\hn$ we define
$$\rvf_k:=\Ra_{\ggzero_0}-\Ra_{\ggtwo_{k}^{k}}\quad , \qquad \cvf_k:=\Ra_{\kgone^{k}}\,.$$
The correspondences $\Ra_k\mapsto \mmat,\ \cvf_k\mapsto \cmat,\ \rvf_k\mapsto \rmat$ imply that the Lie algebra generated by $\Ra_k,\, \rvf_k,\, \cvf_k$ is isomorphic to $\mathfrak{sl}_2(\C)$. Hence there are $\hn$ copies of $\sl2$ in AMSY-Lie algebra $\amsy$ attached to the non-rigid compact \cy threefolds on $\C$.


\def\cprime{$'$} \def\cprime{$'$} \def\cprime{$'$}





\end{document}